\documentclass[11pt]{amsart}
\usepackage{psfrag}

\usepackage{multicol}
\usepackage{array,amscd,amsmath,amsthm}
\usepackage{amssymb}
\usepackage{graphicx}
\usepackage[all]{xy}
\usepackage{multicol}
\allowdisplaybreaks

\makeatletter

\newenvironment{figurehere}
  {\def\@captype{figure}}
  {}
\makeatother

\newtheorem{theorem}{Theorem}[section]
\newtheorem{proposition}[theorem]{Proposition}
\newtheorem{lemma}[theorem]{Lemma}
\newtheorem{corollary}[theorem]{Corollary}

\newtheorem*{proposition*}{Proposition}
\newtheorem*{theorem*}{Theorem}

\theoremstyle{remark}
\newtheorem{definition}[theorem]{Definition}

\newtheorem{remark}[theorem]{Remark}
\newtheorem*{definition*}{Definition}
\newtheorem*{notation}{Notation}

\newcommand{\C}{\mathbb{C}}
\newcommand{\Z}{\mathbb{Z}}

\newcommand{\R}{\mathbb{R}}

\def\extended#1{\tilde{#1}}

\def\Ii{\mathrm{Im}}

\def\M{\mathcal{M}}
\def\eM{\extended{\M}}
\def\Mbar{\overline{\M}}
\def\Cc{\mathcal{C}}

\def\H{\mathbb{H}}
\def\Teich{\mathcal{T}}
\def\eTeich{\extended{\Teich}}

\def\pa{\partial}
\def\Pr{\mathbb{P}}
\def\aut{\mathrm{Aut}}
\def\iso{\mathrm{Iso}}

\def\rar{\rightarrow}
\def\arr#1#2{\stackrel{#1}{#2}}
\def\hra{\hookrightarrow}
\def\a{\alpha}
\def\Qq{\mathcal{Q}}
\def\Af{\mathfrak{A}}

\def\b{\beta}
\def\ua{\mathbf {\a}}

\def\Ao{\Af^\circ}
\def\A8{\Af_{\infty}}

\def\l{\lambda}

\def\i{\iota}
\def\g{\gamma}
\def\G{\Gamma}

\def\d{\delta}
\def\lra{\longrightarrow}

\def\Ll{\mathcal{L}}

\def\s{\sigma}
\def\Si{\Sigma}
\def\e{\varepsilon}

\def\ul#1{\underline{#1}}
\def\ol#1{\overline{#1}}
\def\wh#1{\widehat{#1}}
\def\dis{\displaystyle}
\def\ti#1{\tilde{#1}}

\def\Sp{\mathrm{Sp}}
\def\Pr{\mathbb{P}}
\def\ora#1{\overrightarrow{#1}}
\def\ola#1{\overleftarrow{#1}}
\def\bold#1{\mbox{\boldmath$#1$}}
\def\ua{\ul{\bold{\a}}}

\def\up{\ul{p}}
\def\os#1#2{\overset{#1}{#2}}
\def\Tw{\mathrm{Tw}}
\def\sst{\scriptstyle}
\def\th{\vartheta}

\begin{document}

\title[Triangulated Riemann surfaces with boundary]
{Triangulated Riemann surfaces with boundary
and the Weil-Petersson Poisson structure}

\author{Gabriele Mondello}

\address{Imperial College of London\\
Department of Mathematics, Huxley Building\\
South Kensington Campus\\
London, SW7 2AZ, UK\\
\phantom{x}}

\email{g.mondello@imperial.ac.uk}


\begin{abstract}
Given a hyperbolic surface with geodesic boundary $S$, the lengths
of a maximal system of disjoint simple geodesic
arcs on $S$ that start and end at $\pa S$
perpendicularly are coordinates on the Teichm\"uller space $\Teich(S)$.  
We express the Weil-Petersson Poisson structure of $\Teich(S)$
in this system of coordinates and we prove that it limits pointwise
to the piecewise-linear Poisson structure defined by Kontsevich
on the arc complex of $S$. At the same time, we obtain
a formula for the first-order variation of the distance between
two closed geodesics under Fenchel-Nielsen deformation.
\end{abstract}


\maketitle

%
\begin{section}*{Introduction}
The Teichm\"uller space $\Teich(S)$ of a
compact oriented surface $S$ with marked points
is endowed with a K\"ahler metric,
first defined by Weil using Petersson's pairing of
modular forms.
By the work of Wolpert (\cite{wolpert:elementary},
\cite{wolpert:fenchel-nielsen} and \cite{wolpert:wp-geometry}),
the Weil-Petersson K\"ahler form $\omega_{WP}$ can be neatly
rewritten using Fenchel-Nielsen coordinates.

Algebraic geometers became interested in Weil-Petersson
volumes of the moduli space of curves $\M(S)=\Teich(S)/\Gamma(S)$
since Wolpert \cite{wolpert:homology} showed
that the class of $\omega_{WP}$ is proportional to
the tautological class $\kappa_1$, previously defined
by Mumford \cite{mumford:towards} in the algebro-geometric
setting and then by Morita \cite{morita:characteristic-bull}
in the topological setting.
The reason for this interest relies on the empirical fact
that many problems in enumerative geometry of algebraic curves
can be reduced to the intersection theory of the
so-called tautological classes (namely, $\psi$ and $\kappa$)
on the moduli space of curves.

A major breakthrough in the 1980s and early 1990s was the discovery
(due to Harer, Mumford, Penner and Thurston)
of a cellularization of the moduli space of punctured
Riemann surfaces,
whose cells are indexed by ribbon graphs (also called fatgraphs),
that is finite graphs together with the datum of a cyclic order
of the half-edges incident at each vertex.
To spell it out better, if $S$ is a compact oriented surface
with distinct marked points $c_1,\dots,c_m\in S$
such that the punctured surface $\dot{S}=S\setminus\{c_1,\dots,c_m\}$
has $\chi(\dot{S})<0$, then
there is a homeomorphism between $\M(S)\times\R^m_+$
and the piecewise-linear space $\M^{comb}(S)$ of metrized
ribbon graphs whose fattening is homotopy equivalent to $\dot{S}$.

By means of this cellularization,
many problems could be attacked
using simplicial methods (for instance, the orbifold Euler
characteristic of $\M(S)$ \cite{harer-zagier:euler}
\cite{penner:euler} and
the virtual homological dimension of the mapping class
group $\G(S)$ \cite{harer:virtual}).
A major success was also Kontsevich's proof
\cite{kontsevich:intersection}
of Witten's conjecture \cite{witten:intersection}, which
says that the generating series of the intersection numbers
of the $\psi$ classes on the compactified moduli spaces
satisfies the KdV hierarchy
of partial differential equations.
One of the key steps in Kontsevich's proof was
to explicitly rewrite the $\psi$ classes on
each cell of $\M^{comb}(S)$ in terms of the affine coordinates,
i.e. the lengths
of the edges of the graph indexing the cell.

A different approach to Witten's conjecture was developed
by Mirzakhani \cite{mirzakhani:virasoro}, by noticing that
the intersection numbers appearing in the generating series
can be better understood as Weil-Petersson volumes of
the moduli space of hyperbolic surfaces with geodesic boundaries
of fixed lengths.
Generalizing a remarkable identity of McShane
\cite{mcshane:simple} involving lengths of simple closed
geodesics,
she was able to unfold the integral over $\M(S)$ to an integral
over the space of couples $(\Si,\g)$, where $\Si$ is a hyperbolic
surface with geodesic boundary homotopy equivalent to $\dot{S}$
and $\g\subset\Si$ is a simple closed geodesic,
and then to relate this space to the moduli spaces of
hyperbolic surfaces homeomorphic to $\Si\setminus\g$.
The recursions she obtains are known as Virasoro equations
and (together with string and dilaton equation) are equivalent
to the KdV hierarchy.

In Mirzakhani's approach, hyperbolic surfaces with boundary
play a key role as the recursion is really built on the
process of cutting a surface along a simple closed geodesic.

Back to the cellularization, there is not just one way to
attach a metrized ribbon graph to a Riemann surface.
A first way is due to Harer-Mumford-Thurston (and is described
in \cite{harer:virtual} and \cite{kontsevich:intersection})
and uses existence and uniqueness of quadratic differentials
with closed horizontal trajectories on $\Si$ and double poles of
prescribed quadratic residues at the punctures (see
\cite{strebel:book}). Another way to rephrase it
is the following:
given a Riemann surface $\Si$ with $m$ marked points
$c_1,\dots,c_m$
and positive numbers $p_1,\dots,p_m$,
there exists a unique way to give a metric $g$
(with conical singularities) to the surface
in its conformal class and to dissect $\Si$ into
pointed polygons $(P_i,c_i)$,
such that each $(\dot{P}_i,g|_{\dot{P}_i})$
is isometric to a semi-infinite flat cylinder
of circumference $p_i$ with $c_i$ at infinity.
The boundaries of the polygons describe a ribbon
graph $G$ embedded in $S$ and
the lengths of the sides of the polygons provide
local affine coordinates for the cells indexed by
the isomorphism type of $G$ (as unmetrized ribbon graph).

Now, we are going to describe
a second way to produce a ribbon graph out of a punctured
surface, which uses hyperbolic geometry
and which is due to Penner \cite{penner:decorated}
and Bowditch-Epstein
\cite{bowditch-epstein:triangulations}
(see also \cite{acgh2:book} for a
detailed explanation).

The uniformization theorem endows every
compact or punctured Riemann surface $\dot{\Si}$
(homotopy equivalent to $\dot{S}$)
with $\chi(\dot{\Si})<0$
with a hyperbolic metric of finite volume, so that its punctures
correspond to cusps.
In this case, a decoration of a punctured
Riemann surface $\dot{\Si}$ is a choice of a horoball
$B_i$ at the $i$-th puncture for every $i$.
If the radii $p_i$ of the horoballs
are sufficiently small, then the $B_i$'s
are all disjoint and we can consider the spine of
the truncated surface
$\Si_{\up}:=\Si\setminus\bigcup_i B_i$, that is the locus of points
whose distance from the boundary $\bigcup_i\pa B_i$
is realized by at least two paths. The wanted ribbon
graph is given by this spine, which is a
one-dimensional CW-complex embedded in the surface
with geodesic edges.
Mimicking what done with quadratic differentials,
one could choose the lengths of the edges of the spine as local
coordinates.
This choice works well for a topological treatment,
but for geometric purposes there are more useful options.

A first interesting system of coordinates on the
Teichm\"uller space of decorated surfaces $\Teich(S)\times\R_+^m$
(defined by Penner) is given by
the lengths $\{\ti{a}_i\}$, where $\{\a_i\}$ is a maximal system
of arcs on $S$ (see Section~\ref{subsec:arcs}), $\ti{a}_i(f:S\rar\Si)=
\ell(f(\a_i)\cap\Si_{\up})$ and $f(\a_i)$ is understood to be
the unique geodesic representative in its homotopy class.
Beside its naturality, the interest for these coordinates
is also due to the following.

\begin{theorem}[\cite{penner:volumes}]
Let $\{\a_i\}$ be a maximal system of arcs on the
pointed surface $S$.
The Weil-Petersson form $\omega$ on $\Teich(S)$ pulls back
on $\Teich(S)\times\R_+^m$ to
\[
-\frac{1}{2}
\sum_{t\in H} (d\ti{a}_i\wedge d\ti{a}_j+d\ti{a}_j\wedge d\ti{a}_k+
d\ti{a}_k\wedge d\ti{a}_i)
\]
where $H$ is the set of ideal triangles in $S\setminus\bigcup_i\a_i$
and $(\a_i,\a_j,\a_k)$
are the arcs bounding $t$, in the cyclic
order compatible with the orientation of $t$.
\end{theorem}

On top-dimensional cells of $\M^{comb}(S)$, the system of
arcs dual to the spine is maximal and so the theorem above
expresses the restriction of the Weil-Petersson
form $\omega$ to maximal cells.
This would be enough to integrate all over $\M^{comb}(S)$ if
we knew how to describe top-dimensional cells in the $\ti{a}_i$'s.

On the other hand, cells can be easily described in
a second remarkable system of coordinates.
Penner's simplicial coordinates associated to
the spine are the lengths of the horocyclic segments
that are projections of edges of the spine.
In these coordinates, cells look like straight simplices
but the lengths of the dual arcs cannot be easily
expressed as functions of their simplicial coordinates.

To compute intersection numbers on a compactification of $\M(S)$,
Kontsevich \cite{kontsevich:intersection} integrates over
simplices of maximal dimension in $\M^{comb}(S)$, even though
he used Harer-Mumford-Thurston's construction to produce
the ribbon graph and the lengths of the edges as local
affine coordinates.

Clearly, these systems of coordinates are very different,
but the integration schemes in \cite{penner:volumes}
and \cite{kontsevich:intersection} for volumes of $\M(S)$
are the same as described above.
The reason of this similarity relies on the following
observation.

Let $S$ be a compact oriented surface with $m$ boundary components.
The moduli space $\M^*(S)$ of hyperbolic surfaces
$\Si$ homeomorphic to $S$ together with a choice of a preferred
point on each component of $\pa \Si$ carries a ``Weil-Petersson''
symplectic structure (see \cite{wolpert:symplectic} and
\cite{goldman:symplectic}). Strictly related to
the boundary length function
$\Ll:\M^*(S)\rar \R^m$ is the moment map $\Ll^2/2$
for the natural $(S^1)^m$-action on $\M^*(S)$ (see
\cite{mirzakhani:virasoro}), whose quotient is $\M(S)$.
The symplectic reductions
are exactly the loci $\M(S)(p_1,\dots,p_m)\subset\M(S)$
of hyperbolic surfaces with boundaries of length
$p_1,\dots,p_m>0$, endowed with the Weil-Petersson
symplectic form.

By general considerations on the symplectic reduction,
one can notice that the (class of the)
``symplectic form'' $\Omega=\sum p_i^2\psi_i$
used by Kontsevich represents the normalized limit
of the Weil-Petersson form $\omega$ on
$\M(S)(p_1,\dots,p_m)$
as $(p_1,\dots,p_m)$ diverges.
Thus, Kontsevich also computed
suitably normalized Weil-Petersson volumes.
On the other hand, decorated surfaces can be thought of
as Riemann surfaces with infinitesimal boundaries.\\

In this paper, we define a natural Poisson structure $\eta$
on the Teichm\"uller space of Riemann surfaces
with boundary using the doubling construction.
Results of Wolpert \cite{wolpert:symplectic}
and Goldman \cite{goldman:symplectic} imply that
the associated bivector field on $\Teich(S)$ has the form
\begin{equation*}\tag{*}\label{eq:wolpert_formula}
\eta_S=-\sum_i \frac{\pa}{\pa\ell_i}\wedge \frac{\pa}{\pa\tau_i}
\end{equation*}
where the sum ranges
over a maximal system of disjoint simple closed curves,
that are not boundary components.
As before,
the symplectic leaves of this Poisson structure are the loci
$\Teich(S)(p_1,\dots,p_m)$ of surfaces with fixed
boundary lengths $p_1,\dots,p_m$, endowed with the Weil-Petersson
symplectic structure (in \cite{goldman:dynamics} it is shown
that this happens more generally for spaces of representations
of $\pi_1(S)$ inside a Lie group).

\begin{remark}
As noted by the referee, if $\Si$ is a compact hyperbolic surface
with no boundary,
Wolpert's formula (\ref{eq:wolpert_formula})
descends from the more basic symplectic duality
$\omega(\pa/\pa\tau_\xi,-)=d\ell_\xi$ and from Fenchel-Nielsen
coordinates.
It is not clear whether the next theorem can descend from
an analogous ``duality'' for hyperbolic surfaces with
boundary.
\end{remark}

Given a hyperbolic surface with geodesic boundary $\Si$,
we can immediately take its spine and so produce a ribbon
graph with no need of decorations.
Clearly, if $\{\a_i\}$ is a maximal system of arcs on $S$,
then the hyperbolic lengths $\{a_i\}$
defined as $a_i(f:S\rar\Si)=\ell(f(\a_i))$
are coordinates
on $\Teich(S)$ (Ushijima \cite{ushijima:decomposition})
and one can check that the difference $a_i-a_j$
limits to the $\ti{a}_i-\ti{a}_j$ for all $i,j$
as the $p_k$'s converge to zero.
More interestingly, Luo \cite{luo:decomposition} showed that
the lengths of the projections of the edges of the spine
to $\pa S$ (which we will call ``widths'') are also
coordinates, which in fact specialize to simplicial coordinates
for infinitesimal boundary lengths (under a suitable normalization).

Our goal is to rewrite the Weil-Petersson Poisson
structure $\eta$ in terms of the $a_i$'s. Our main result is
the following.

\begin{theorem*}
Let $S$ be a compact oriented surface with $m$ boundary
components $\Cc$
and let $\ua=\{\a_1,\dots,\a_{6g-6+3m}\}$ be a maximal system of
arcs of $S$.
The Weil-Petersson Poisson structure on $\Teich(S)$
at $[f:S\rar\Si]$ can be written as
\[
\eta_S=\frac{1}{4}\sum_{C\in\Cc} \sum_{\substack{y_i\in f(\a_i\cap C) \\
y_j\in f(\a_j\cap C) }}
\frac{\sinh(p_C/2-d_C(y_i,y_j))}{\sinh(p_C/2)}\,
\frac{\pa}{\pa a_i}\wedge\frac{\pa}{\pa a_j}
\]
where $p_C=\ell_C$ and $d_C(y_i,y_j)$ is
the length of geodesic arc running from $y_i$ to $y_j$ along
$f(C)$ in the positive direction.
\end{theorem*}

The proof of the theorem above relies on the formula $(*)$
and on the understanding of how the distance between two geodesics
in a surface $R$ without boundary
(we will then take $R$ to be the double of $\Si$)
varies at first order, when we perform a Fenchel-Nielsen deformation.
Let us recall that the (right)
Fenchel-Nielsen deformation along a simple closed
geodesic $\xi$ of $R$ is obtained by cutting $R$ along $\xi$,
letting the left component slide forward of $t$ and then
reglueing the left with the right part.
According to Thurston, it is called ``right'' because
one jumps to the right when one crosses the fault line.
We call $\pa/\pa\tau_{\xi}$
its associated vector field on $\Teich(R)$.

The following result (which we state in a simplified
version, for brevity) might be interesting on its own:
a more complete statement (Theorem~\ref{thm:derivative2})
can be found in Section~\ref{subsec:general}.
It should be compared to Theorem~3.4 of \cite{wolpert:symplectic}.

\begin{theorem*}
Let $R$ be a hyperbolic surface without boundary
and $\delta\subset R$ a simple geodesic
arc running from $y_1\in\g_1$ to $y_2\in\g_2$ that realizes the
distance between the geodesics $\g_1$ and $\g_2$ in its homotopy class.
{\rm Assume that $\xi$ does not intersect $\delta$ and that
no portion of $\xi$ is homotopic to $\delta$.}
Then
\[
\frac{\pa}{\pa\tau_{\xi}}(h) = c_1+c_2
\]
where $\dis c_i=\sum_{x_i\in\xi\cap\g_i} c_i(x_i)$
and
\[
c_i(x_i)=
\begin{cases}
\dis \mathrm{sgn}(d(y_i,x_i))\frac{\exp[-|d(y_i,x_i)|]}{2}\sin(\nu_{x_i}) &
\text{if $\g_i$ is open} \\
\phantom{x}\\
\dis\frac{\sinh(p_i/2-d(y_i,x_i))}{2\sinh(p_i/2)} \sin(\nu_{x_i})
& \text{if $\g_i$ is closed}
\end{cases}
\]
where $h$ is the length of $\delta$, $p_i$ is the length of $\g_i$
(if $\g_i$ is closed),
$\nu_i$ is the angle of intersection at $x_i$ between $\xi$ and $\g_i$.
If $\g_i$ is open, then
$d(y_i,x_i)$ is the distance with sign between $y_i$ and $x_i$
along $\g_i$; if $\g_i$ is closed, then we set $d(y_i,x_i)\in(0,p_i)$.
\end{theorem*}

As suggested by the referee, it seems that the same methods can be
employed to obtain a formula for the second twist derivative.
Though reasonable,
the upshot looks quite complicated, so we will not pursue this
calculation here.

As an easy corollary of our main theorem, we obtain that Kontsevich's
piecewise-linear form $\Omega$ on $\M^{comb}(S)$ that represents
the class $\sum_i p_i^2\psi_i$ is the {\it pointwise limit}
(under a suitable normalization) of twice the Weil-Petersson
form $2\omega_{(p_1,\dots,p_m)}$ on $\M(S)(p_1,\dots,p_m)$
as $(p_1,\dots,p_m)\lra+\infty$.\\

Quite recently, Carfora-Dappiaggi-Gili \cite{CDG:twistor} have
found a different procedure to relate decorated hyperbolic surfaces,
``decorated'' flat surfaces with conical points and hyperbolic
surfaces with geodesic boundary components. It would be interesting
to understand how it relates to the constructions that
we employ here.
%
%
\begin{subsection}*{Plan of the paper}
Section~\ref{sec:boundary} deals with preliminary
results on Riemann surfaces $S$ with boundary,
the construction of the real double $dS$ and intrinsic metrics.
We recall the definition of Teichm\"uller space $\Teich(S)$
and Weil-Petersson form, and we establish a link between
the Poisson structure of $\Teich(dS)$ and that of $\Teich(S)$.
We also recall the definition of arc complex of a surface with boundary.

In Section~\ref{sec:triangulations}, we describe the construction
of the spine and we illustrate the results of Ushijima and Luo,
who define two different system of coordinates on $\Teich(S)$
using triangulations of $S$, and we show how to decompose $\Teich(S)$
into ideal cells in a $\Gamma(S)$-equivariant way.
We compare their theorems to previous results of Penner, who proved
the analogous statements 
in the case of decorated Riemann surfaces, and we
show that Ushijima and Luo's coordinates specialize to Penner's ones.

In Section~\ref{sec:FN}, we review the Fenchel-Nielsen deformation and
we use techniques of Wolpert to compute
the first-order variation of the distance
between two geodesics.

In Section~\ref{sec:poisson}, we establish our main result and write
the Weil-Petersson Poisson structure in the coordinates $\{a_i=\ell_{\a_i}\}$
for every maximal system of arcs $\{\a_i\}$ on $S$, using
Wolpert's formula and the result of Section~\ref{sec:FN}.
As a corollary, we deduce that Kontsevich's PL representative for $\Omega$
is the pointwise limit of the Weil-Petersson form, when the boundary
components become infinitely large.
We also check that our result agrees with Penner's computations
for decorated surfaces.

Appendix~\ref{app:hyperbolic} collects a few formulae of
elementary hyperbolic trigonometry, that are used in the rest
of the paper.
\end{subsection}
%
%
\begin{subsection}*{Acknowledgements}
I would like to thank Enrico Arbarello,
William Goldman, Curtis McMullen, Tomasz Mrowka, Robert Penner
and Scott Wolpert for useful discussions.
In particular, I thank Francesco Bonsante for suggesting to me
to use Wolpert's techniques in the proof of Theorem~\ref{thm:derivative}
and for many fruitful conversations.
Finally, I am grateful to an anonymous referee 
for carefully reading this paper and for helpful remarks
and comments.
\end{subsection}
%
%
\end{section}
%
%
\begin{section}{Riemann surfaces with boundary}\label{sec:boundary}
\begin{subsection}{Double of a surface with boundary}
A compact {\it surface with nodes} and boundary
is a compact Hausdorff topological
space $R$ with countable basis that is locally homeomorphic
either to $\C$, or to $\{z\in\C\,|\, \mathrm{Im}(z)\geq 0\}$
or to $\{(z,w)\in\C^2\,|\,zw=0\}$.
Points of $R$ that have a neighbourhood of the first type
are said {\it smooth}; in the second case,
points on the real line are said to belong to the
{\it
boundary}; in the third case,
the point $\{z=w=0\}$ is called
a {\it node}.
We will always assume that $\pa R$ is homeomorphic to a disjoint union of $b$
copies of $S^1$, that $R$ is connected (unless
differently specified) and that $R$ is always endowed
with the unique differentiable structure away from the nodes.

The (arithmetic) {\it genus} of such a connected surface is
$g=1+(\nu-\chi-b)/2$, where $\chi$ is the Euler characteristic,
$b$ is the number of boundary components and and $\nu$ is the number of nodes.\\

Consider a compact oriented surface $\Si$ of genus $g$ 
with boundary circles $C_1,\dots,C_n$ (endowed with the
orientation induced by $\Si$)
and let $c_1,\dots,c_m\in \Si$ be distinct smooth
{\it marked} points.
We will also write $\dot{\Si}$ for the {\it punctured} surface
$\Si\setminus\{c_1,\dots,c_m\}$.

Call $\Si'$ the oriented surface obtained from $\Si$ switching
the orientation and similarly denote by $C'_1,\dots,C'_n$ its
boundary components and $c'_1,\dots,c'_m$ its marked points.
In general, for every point $x\in \Si$
call $x'$ the corresponding point in $\Si'$.

The {\it double} of $\Si$ is the compact oriented
surface $d\Si$ of arithmetic
genus $2g+(n+m)-1$ without boundary obtained from $\Si \sqcup \Si'$
identifying $x\sim x'$ for every $x\in\pa \Si\cup\{c_1,\dots,c_m\}$.
Clearly, $d\Si$ is connected if $n+m>0$ and it has nodes if $m>0$.
Call $\i:\Si\hra d\Si$ and $\i':\Si'\hra d\Si$ the natural inclusions.

The surface $d\Si$ has a natural orientation-reversing
involution $\s$ which exchanges $\i(\Si)$ with $\i'(\Si')$
and fixes $\i(\{c_1,\dots,c_m\})$ and $\i(\pa \Si)$ pointwise.

Suppose now that $d\Si$ has a complex-analytic structure $J$,
meaning that the nodes of $d\Si$ have a neighbourhood
biholomorphic to $\{(z,w)\in\C^2\,|\,zw=0,\, |z|<\e,\, |w|<\e\}$.
We say that $J$ is {\it compatible} with the involution $\s$
if the homeomorphism
$\s:d\Si\rar d\Si$ is anti-holomorphic,
or in other words $\s^* J=-J$.
This implies that $\i(\pa \Si)$ is a totally real submanifold
(and the $\i(c_j)$'s are real points) of $(d\Si,J)$.
Conversely, an atlas of charts on $\Si$,
which are holomorphic on $\Si^\circ:=\Si\setminus\pa \Si$ and
map the boundary of $\Si$ to $\mathbb{R}\subset\mathbb{C}$
and the marked points to $0\in\mathbb{C}$, is the restriction through $\i$
of a complex structure on $d\Si$ compatible with $\s$.
In this case, we will say that $(\Si,\i^* J)$ is a
{\it Riemann surface with boundary} and we will denote it
just by $\Si$ when the complex structure is understood.

A {\it morphism} between Riemann surfaces with boundary is a continuous
application $\Si_1\rar \Si_2$ that maps
$\pa \Si_1$ to $\pa \Si_2$ in a real-analytic way
and that restricts to a holomorphic map
$\Si_1^\circ\rar \Si_2^\circ$ on the interior,
preserving the marked points.
Equivalently, it is the restriction of a holomorphic map
between the doubles $d\Si_1\rar d\Si_2$ that
commutes with the $\s$-involutions.
\end{subsection}
\begin{subsection}{Metrics on a Riemann surface with
boundary}\label{subsec:metrics}
We can associate two natural metrics to a smooth Riemann
surface $\Si$ of genus $g$
with $n$ boundary components and $m$ marked points,
if $\chi(\dot{\Sigma})<0$
(which actually coincide if $n=0$).

This is a consequence of
the uniformization theorem, which says that the universal cover
of a Riemann surface is biholomorphic either to the Riemann sphere
$\C\Pr^1$, or to the complex plane $\C$
or to the upper half-plane $\H$.
The crucial fact is that $\H$
has a complete Hermitean metrics of constant curvature
(the hyperbolic metric
$y^{-2}{dz\,d\bar{z}}$)
that is preserved by all analytic automorphisms.
Similarly, the Fubini-Study metric on $\C\Pr^1$
is complete, of constant positive curvature and
invariant under automorphisms that preserve the real line $\R\Pr^1$.
The flat metric $dz\,d\bar{z}$ on $\mathbb{C}$ is canonically
defined up to rescaling.
We call these metrics {\it standard}.

Back to $\Si$, the first natural metric is defined
considering $\dot{\Si}^\circ$ as an open Riemann surface:
excluding the case $\Si\cong\C\Pr^1$,
the universal cover $\ti{\dot{\Si}}^\circ$ is
isomorphic to $\C$ if $(g,n,m)=(0,0,1),(0,0,2),(1,0,0)$, and to
$\H$ otherwise.
In the last case,
the covering map $\ti{\dot{\Si}}^\circ\rar \dot{\Si}^\circ$ determines
a {\it holonomy} representation $\rho:\pi_1(\dot{\Si}^\circ)\hra
\aut(\ti{\dot{\Si}}^\circ)=
\iso(\ti{\dot{\Si}}^\circ)$, uniquely defined up to inner automorphisms
of $\aut(\ti{\dot{\Si}}^\circ)$. The standard Hermitean metric on $\ti{\dot{\Si}}^\circ$
descends to a complete
Hermitean metric on $\ti{\dot{\Si}}^\circ/\Ii(\rho)\cong \dot{\Si}^\circ$.

The second natural metric on $\Si$
is obtained by restricting the standard metric on the double $d\Si$
via the inclusion $\i$.
The universal cover of (each connected component of) the smooth locus
$d\Si_{sm}$ of $d\Si$ is isomorphic to $\H$
if $(g,n+m)\neq (0,0),(0,1),(0,2),(1,0)$.
In these cases,
$\dot{\Si}$ inherits a complete Hermitean
metric of curvature $-1$ with totally geodesic boundary $\pa \Si$,
which is called the {\it intrinsic metric}. Under a suitable
normalization (which fixes the curvature or the area), it is uniquely
determined by the isomorphism class of $\dot{\Si}$.
Both metrics
acquire cusps at the marked points.
\end{subsection}
\begin{subsection}{The extended Teichm\"uller space}
Fix a compact oriented smooth
surface $S$ of genus $g$ with boundary
components $C_1,\dots,C_n$ and let $\Si$
be a smooth Riemann surface, possibly
with boundary and marked points.
An {\it $S$-marking} of $\Si$ is a
smooth map $f:S\lra\Si$ that may contract boundary components to
marked points
such that $f_{int}:S^\circ\lra\dot{\Si}^\circ$ is an
orientation-preserving diffeomorphism.

The {\it extended Teichm\"uller space} of $S$ is the space $\eTeich(S)$
the space of equivalence classes of $S$-marked
Riemann surfaces
\[
\eTeich(S):=\{f:S\lra \Si \,|\, \Si\ \text{Riemann surface}\}/\sim
\]
where $f:S\lra \Si$ is an {\it $S$-marking}
and the equivalence relation $\sim$ identifies
$f$ and $f':S\arr{\sim}{\lra} \Si'$
if and only if there exists an isomorphisms of
Riemann surfaces $h: \Si\arr{\sim}{\lra} \Si'$ such that
$(f'_{int})^{-1}\circ h \circ f_{int}$ is isotopic to the identity.
The Teichm\"uller space $\Teich(S)\subset\eTeich(S)$ is the locus
of those class of 
markings $f:S\rar \Si$ which do not shrink any boundary component
to a point.

There are several ways to put a topology on $\eTeich(S)$.
For instance, we have seen in Section~\ref{subsec:metrics} that
a complex structure on $\Si$ determines and is determined
by a complete hyperbolic metric on $\dot{\Si}$ with totally
geodesic boundary.
The universal cover of $\dot{\Si}^\circ$ has a developing map
into $\H$ and so a holonomy map $\pi_1(\dot{\Si}^\circ)\rar
\mathrm{PSL}_2\R$ is induced. Pulling it back through the marking
$f_{int}:S^\circ\rar\dot{\Si}^\circ$, we get a global injection
(originally due to Fricke)
\[
\xymatrix{
\eTeich(S) \ar@{^(->}[rr] &&
\mathrm{Hom}(\pi_1(S),\mathrm{PSL}_2 \R)/\mathrm{PSL}_2 \R
}
\]
which is independent of the choices made: thus, we
can endow $\eTeich(S)$ with the subspace topology.

Hence, the Teichm\"uller space of $S$
can be thought of as the space of complete hyperbolic
metrics on $S$ with totally geodesic boundary (up to isotopy).
Points in $\eTeich(S)\setminus\Teich(S)$ correspond to $S$-marked
hyperbolic surfaces in which some boundary components of $S$
are collapsed to cusps of $\Si$.

Thus, we have a natural {\it boundary-lengths map}
\[
\xymatrix@R=0pt{
\hspace{-1cm}\Ll\ :\quad
\eTeich(S) \ar[r] & \R_{\geq 0}^n \\
[f:S\rar (\Si,g)] \ar@{|->}[r] &
\dis \left( \ell_{C_1}(f^*g),\dots,\ell_{C_n}(f^*g)
\right)
}
\]
If we call $\eTeich(S)(p_1,\dots,p_n)$ the submanifold
$\Ll^{-1}(p_1,\dots,p_n)$, then $\Teich(S)=\eTeich(S)(\R_+^n)$.

Define the {\it mapping class group}
$\Gamma(S)$ as $\pi_0\mathrm{Diff}_+(S)$, that is the group
of orientation-preserving diffeomorphisms of $S$ that send each
boundary component to itself, up to isotopy.

The group $\Gamma(S)$ acts properly and discontinuously on $\eTeich(S)$:
its quotient $\eM(S):=\eTeich(S)/\Gamma(S)$ is
the {\it extended moduli space of Riemann surfaces with boundary}.
The moduli space itself $\M(S)=\Teich(S)/\Gamma(S)\subset \eM(S)$
is naturally an orbifold.
\end{subsection}
%
\begin{subsection}{Deformation theory of Riemann surfaces with boundary}
Let $S$ be a smooth
compact Riemann surface with boundary and $\chi(S)<0$,
and let $[f:S\rar\Si]\in\eTeich(S)$.
We want to understand the deformations of
$\Si$ as a Riemann surface
with boundary (and possibly cusps).
We refer to \cite{Deligne-Mumford:irreducibility},
\cite{Bers:degenerating} and \cite{bers:deformations}
for a more detailed treatment of
the case of surfaces with nodes.

A first way to approach the problem is to pass to its double $d\Si$.
Suppose first that $\Si$ has no cusps and so $d\Si$ is smooth.

The space of first-order deformations of complex structure
on the surfaces $d\Si$ can be identified
to the complex vector
space $\mathcal{H}(d\Si)$ of harmonic Beltrami differentials.
If $g$ is the hyperbolic metric on $d\Si$ and
$\Qq(d\Si)$ is the space of holomorphic quadratic differentials
on $d\Si$ (i.e. holomorphic sections of $(T^*_{d\Si})^{\otimes 2}$),
then the elements of $\mathcal{H}(d\Si)$ are
$(0,1)$-forms $\mu$ with values in the tangent bundle of $\Si$
which are harmonic with respect to $g$,
so that they can be written as $\mu=\ol{\varphi}/g$,
for a suitable $\varphi\in\Qq(d\Si)$. Thus, $\mathcal{H}(d\Si)$
can be identified to the dual of $\Qq(d\Si)$.

To construct complex-analytic charts of $\Teich(dS)$
(or of $\M(dS)$), one can use a method
mostly due to Grothendieck and whose details can be found
in \cite{arbarello-cornalba:kuranishi} and in \cite{acgh2:book}.
It relies on the fact that smooth
compact Riemann surfaces,
with negative Euler characteristic, can be embedded through
the tricanonical linear system
in a complex projective space. Thus, all holomorphic
families of such curves can be pulled back from
a smooth open subset $\mathcal{V}$ of a Hilbert scheme and
a semi-universal deformation $\mathcal{D}(d\Si)$
of $d\Si$ (which means that the Kodaira-Spencer map is an
isomorphism at every point of $\mathcal{D}(d\Si)$)
can be obtained just taking
a suitable slice of $\mathcal{V}$.
After restricting the family over a ball,
$\mathcal{D}(d\Si)$ gives an complex-analytic
orbifold chart for $[d\Si] \in \M(dS)$ and 
(choosing a smooth trivialization of the family)
an honest chart for a neighbourhood of any $[f:dS\rar d\Si]
\in \Teich(dS)$.

Because $\s$ acts on $\mathcal{D}(d\Si)$
as an antiholomorphic involution, then
the first-order deformations of complex structure
on $dS$ compatible with the $\s$-involution
are parametrized by the real subspace
$\mathcal{H}(d\Si)^{\s}$, dual to $\Qq(d\Si)^{\s}$
(which can be identified to the real vector space of
holomorphic quadratic differentials on $\Si$, whose
restriction to $\pa\Si$ is real).

If $\Si$ has $k$ cusps, then $d\Si$ has $k$ nodes $\nu_1,\dots,\nu_k$
and the semi-universal deformation $\mathcal{D}(d\Si)$
(and so the orbifold chart for $\Mbar(dS)$ around $[d\Si]$) can
still be constructed slicing the Hilbert scheme of curves
embedded using the third power of their dualizing line bundle.

If $\mathcal{S}\rar \mathcal{D}(d\Si)$ is the tautological family,
one can find an open subset $U_i\subset\mathcal{S}$
with local analytic coordinates $z_i,w_i$
such that the deformation of the node $\nu_i$ looks like $\{z_i w_i=t_i\}$,
where $\{t_1,\dots,t_k,s_1,\dots,s_N\}$ are a system of coordinates
at $[d\Si]$ on $\mathcal{D}(d\Si)$.

The smooth divisor $\mathcal{N}_i=\{t_i=0\}
\subset\mathcal{D}(d\Si)$ parametrizes those deformations of $d\Si$
in which the node $\nu_i$ survives. Call $\mathcal{N}=\bigcap_{i=1}^k
\mathcal{N}_i$.

As a consequence, the space of first-order deformations of $d\Si$
is given by
\[
\xymatrix@R=0.25in{
0 \ar[r] & T_{d\Si}\mathcal{N} \ar[r] &
T_{d\Si}\mathcal{D}(d\Si) \ar[r] & N_{\mathcal{N}/
\mathcal{D}(d\Si)} \ar[r] & 0 \\
& \mathcal{H}(d\Si) \ar[u]_\cong && \C^k \ar[u]_\cong
}
\]
where in this case $\mathcal{H}(d\Si)$ is the space of harmonic
Beltrami differentials on $d\Si$ that vanish at the nodes
and $\C^k$ is spanned by the $\pa/\pa t_i$'s.

Consequently, the space of first-order deformations
of $\Si$ is given by
\[
\xymatrix{
0 \ar[r] & \mathcal{H}(d\Si)^{\s} \ar[r] &
T_{\Si}\mathcal{D}(\Si) \ar[r] & \R^k \ar[r] & 0
}
\]
where $\R^k=(\C^k)^{\s}$. However, only the directions
that project to $(\R_{\geq 0})^k\subset\R^k$
(corresponding to $t_1,\dots,t_k\geq 0$)
belong to the tangent cone. In fact, being interested
in the deformations of $d\Si$ that preserve the symmetry $\s$,
we can choose $w_i=\ol{z}_i$, and so $t_i=z_i w_i=|z_i|^2\geq 0$.

From a different perspective (using harmonic maps),
it follows from \cite{wolf:degeneration} that
the tangent cone to $\eTeich(S)$ at $[f:S\rar\Si]$
can be parametrized by the space $\Qq(\Si)$ of
quadratic differentials, which are holomorphic
on $\dot{\Si}$,
real along the boundary components and that look like
$(a_{-2}^2 z^{-2}+a_{-1}z^{-1}+\dots)dz^2$ at the cusps, with $a_{-2}\leq 0$.

Both approaches show that $\eTeich(S)$ can be made into a
real-analytic smooth variety with corners.

From a global point of view,
the Teichm\"uller space $\Teich(S)$ has a natural
embedding $D:\Teich(S) \hra \Teich(dS)$ 
onto the real-analytic submanifold
of $dS$-marked Riemann surfaces that carry an anti-holomorphic
involution isotopic to $\s$.

\begin{remark}
The inclusion above can be extended to
an embedding $\ol{D}:\ol{\Teich}(S)\hra\ol{\Teich}(dS)$, where
$\ol{\Teich}$ (which contains $\eTeich$)
is the Deligne-Mumford bordification of $\Teich$.
We will not deal with $\ol{\Teich}$: for further details,
see also \cite{looijenga:cellular} or \cite{acgh2:book}.
\end{remark}
\end{subsection}
%
%
\begin{subsection}{The Weil-Petersson form}
Let $S$ be a smooth
compact Riemann surface with boundary with $\chi=\chi(S)<0$ and
consider the universal family $\pi:\mathcal{S}\lra\eTeich(S)$
over the Teichm\"uller space of $S$.
The fibers of $\pi$ are $S$-marked surfaces endowed
with a metric of constant negative curvature $-1$, that is
the vertical tangent bundle $T_\pi$ over $\mathcal{S}$ is
endowed with a Hermitean metric $g$.
\begin{definition}
The {\it Weil-Petersson bivector field} on
$\eTeich(S)$
at $[f:S\rar(\Si,g)]$ is given by
\[
\eta_S(\varphi,\psi):=\mathrm{Im}
\int_{\Si} \frac{\varphi \bar{\psi}}{g}
\]
for every $\varphi,\psi\in \Qq(\Si)\cong T^*_{[f]}\eTeich(S)$.
\end{definition}

Clearly, we can also easily
define the {\it Weil-Petersson $2$-form $\omega_S$} on
$\Teich(S)$ at $[f:S\rar(\Si,g)]$ (i.e. where $\Si$
acquires no cusps) as
\[
\omega_S(\mu,\nu):=\mathrm{Im}
\int_{\Si} \mu \bar{\nu}\cdot g
\]
for $\mu,\nu\in \mathcal{H}(\Si)\cong T_{[f]}\Teich(S)$.

\begin{remark}
The divergence occurring when $\Si$ acquires cusps,
that is when $d\Si$ acquires nodes, was first shown by
Masur (\cite{masur:wp}) using local
coordinates due to Earle and Marden.
As one can notice below, the Weil-Petersson
form is smooth in Fenchel-Nielsen coordinates.
Thus, the differentiable structure
of $\overline{\mathcal{M}}_{g,n}$
underlying the complex-analytic one is
different from the Fenchel-Nielsen differentiable structure; a phenomenon
that was investigated more deeply by Wolpert (\cite{wolpert:wp-geometry}).
\end{remark}

There is another way to describe the Weil-Petersson form on $\eTeich(S)$.
A pair of pants decomposition of $S$ determines Fenchel-Nielsen
coordinates $\ell_1,\dots,\ell_{3g-3+n}\in\R_+$,
$\tau_1,\dots,\tau_{3g-3+n}\in\R$
and also $p_i=\ell_{C_i}\geq 0$ for every boundary component $C_i$
of $S$.

\begin{theorem}[\cite{wolpert:symplectic},
\cite{goldman:symplectic}]\label{thm:wolpert}
The Weil-Petersson $2$-form can be written as
\[
\omega_S=\sum_{i=1}^{3g-3+n}d\ell_i\wedge d\tau_i
\]
on $\eTeich(S)$, with respect to any pair of pants decomposition.
\end{theorem}

\begin{remark}
Literally, Wolpert proved Theorem~\ref{thm:wolpert}
for closed Riemann surfaces, but an inspection of his
paper \cite{wolpert:symplectic} shows that the
statement holds also for Riemann surfaces with boundary.

In \cite{goldman:symplectic}, Goldman defines the Weil-Petersson
symplectic form on the representation variety
of a closed surface. The same definition and treatment
can be extended
to the representation variety of nonclosed surfaces with
or without prescribed holonomy along the boundary components
(see, for instance, \cite{goldman:dynamics}).
\end{remark}

As a consequence, if $S$ is a closed surface,
then $(\Teich(S),\omega_S)$ is a symplectic manifold.
If $S$ has $n$ boundary components,
then $\omega_S$ is degenerate on $\eTeich(S)$, but
$(\eTeich(S)(p_1,\dots,p_n),\omega_S)$ is a symplectic manifold
for all $p_1,\dots,p_n\geq 0$.
\end{subsection}
%
%
\begin{subsection}{Double of a Riemann surface and Weil-Petersson
Poisson structure}\label{sec:2vector}
Consider a smooth compact hyperbolic
Riemann surface $S$ of genus $g$ with boundary components
$C_1,\dots,C_n$ and let $dS$ be its double.

It follows directly from the definition that
the embedding $D:\Teich(S)\hra\Teich(dS)$ induced
by the doubling construction is Lagrangian.
Hence, we relate the Weil-Petersson structures on
$\Teich(S)$ and $\Teich(dS)$ in a different way.

There is a natural map $\pi_\i:\Teich(dS)\lra\Teich(S)$
induced by the inclusion $\i:S\hra dS$
that associates to $[f:dS\arr{\sim}{\lra} (R,g)]$
the $S$-marked hyperbolic subsurface of $R$ with geodesic
boundary isotopic to $f(\i(S))$.

Call $T\Teich(dS)\Big|_{\Teich(S)}$ the restriction
of the tangent bundle of $\Teich(dS)$ through $D$.

\begin{definition}
Set $\displaystyle
\hat{\eta}_S:=(\pi_\i)_*\left(\eta_{dS}\Big|_{\Teich(S)}\right)$, where
$\eta_{dS}$ is the Weil-Petersson bivector field on $\Teich(dS)$ and
$(\pi_\i)_*:T\Teich(dS)\Big|_{\Teich(S)}\lra T\Teich(S)$.
\end{definition}

\begin{proposition}\label{prop:2vector}
The bivector field $\hat{\eta}_S$ coincides with $\eta_S$
on $\Teich(S)$ and we can extend $\hat{\eta}_S$ to $\eTeich(S)$
by setting it equal to $\eta_S$, so that they
define a Poisson structure on
$\eTeich(S)$, whose symplectic leaves are the fibers
of $\Ll:\eTeich(S)\rar\R_{\geq 0}^n$.
\end{proposition}

\begin{proof}
The bivector $\hat{\eta}_S$ defines a Poisson structure
on $\Teich(S)$ because it is obtained pushing $\eta_{dS}$ forward
and $\eta_{dS}$ defined a Poisson structure on $\Teich(dS)$.
The equality $\eta_S=\hat{\eta}_{S}$ follows from Wolpert's
work \cite{wolpert:symplectic}.

To verify this second claim, pick a pair of pants decomposition for $S$.
On $\Teich(dS)$
we have Fenchel-Nielsen coordinates $\ell_i,\tau_i,
\ell'_i,\tau'_i$ for $1\leq i\leq 3g-3+n$
plus $(p_j,\hat{\tau}_j)$,
where $p_j=\ell_{\i(C_j)}$ and $\hat{\tau}_j$ is the
twist parameter of $\i(C_j)$.
By Theorem~\ref{thm:wolpert} we have
\[
\eta_{dS}=-\sum_{i=1}^{3g-3+n}
\left(
\frac{\pa}{\pa \ell_i}\wedge
\frac{\pa}{\pa \tau_i}-
\frac{\pa}{\pa \ell'_i}
\wedge
\frac{\pa}{\pa \tau'_i}
\right)-
\sum_{j} \frac{\pa}{\pa p_j}\wedge
\frac{\pa}{\pa \hat{\tau}_j}
\]
because switching orientation changes the sign
of the twist. Hence
\[
(\pi_\i)_*\left(\eta_{dS}\Big|_{\Teich(S)} \right)=
-\sum_{i=1}^{3g-3+n}
\frac{\pa}{\pa \ell_i}\wedge
\frac{\pa}{\pa \tau_i}
\]
which is vertical with respect to $\Ll$
and whose restriction to each fiber of $\Ll$
is dual to the Weil-Petersson form
according to Theorem~\ref{thm:wolpert}.
\end{proof}
\end{subsection}
%
%
\begin{subsection}{The complex of arcs}\label{subsec:arcs}
Let $S$ be a smooth compact Riemann surface with boundary
components $C_1,\dots,C_n$ and marked points $c_1,\dots,c_m$.
Assume $n+m>0$.

An {\it arc} on $S$ is an embedded unoriented path with endpoints in
$\pa S\cup\{c_1,\dots,c_m\}$,
which is homotopically nontrivial relatively to $\pa S\cup\{c_1,\dots,c_m\}$.
A {\it $k$-system of arcs} is a set of $k$ arcs that are allowed
to intersect only at the marked points of $S$, and which are
pairwise nonhomotopic (relatively to $\pa S\cup\{c_1,\dots,c_m\}$).
We will always consider arcs and systems of arcs up to isotopy
of systems of arcs.

\begin{definition}
The {\it complex of arcs} $\Af(S)$ of $S$ is the simplicial complex,
whose $k$-simplices are $(k+1)$-systems of arcs on $S$.
Maximal simplices of $\Af(S)$ are called {\it triangulations}.
\end{definition}

The complex of arcs was introduced by Harer in \cite{harer:virtual}
(see also \cite{looijenga:cellular}).

A systems of arcs $\ua=\{\a_0,\dots,\a_k\}\in A(S)$ {\it fills} if
$S\setminus\bigcup_i \a_i$ is a disjoint union of discs;
$\ua$ {\it quasi-fills} if $S\setminus\bigcup_i \a_i$ is a disjoint
union of discs, pointed discs
and annuli that retract onto a boundary component.
Call $\Ao(S)\subset \Af(S)$ the subset of systems that quasi-fill
and $\A8(S)\subset \Af(S)$ the subset of those that do not:
$\A8(S)$ is a subcomplex of $\Af(S)$.
Write $|\Ao(S)|$ for $|\Af(S)|\setminus|\A8(S)|$, which is open
and dense inside $|\Af(S)|$.

Also, define $|\Ao(S)|_{\R}:=|\Ao(S)|\times\R_+$.
The space of ribbon graphs $\M^{comb}(S)$ mentioned in the introduction
is homeomorphic to $|\Ao(S)|_{\R}/\Gamma(S)$: given a system of arcs
$\ua$ that quasi-fills, we can construct a ribbon graph embedded in $S$
by drawing edges transversely to the arcs of $\ua$.
In this duality, the weight of an arc corresponds to the length of
its dual edge in the ribbon graph.

%
%
%
%
\begin{remark}
If $\Si$ is a hyperbolic surface, by an {\it arc} $\a$ on $\Si$ we
will usually mean the unique geodesic arc in the isotopy class
of $\a$ that meets $\pa\Si$ perpendicularly, unless differently
specified.
\end{remark}
\end{subsection}
\end{section}
%
%
\begin{section}{Triangulations and spines}\label{sec:triangulations}
Let $S$ be a compact
hyperbolic Riemann surface with nonempty boundary.
For every arc $\a$ on $S$, we have the length function
$\ell_{\a}:\eTeich(S)\rar (0,+\infty]$
that associates to $[f:S\rar (\Si,g)]$
the length of the arc $f(\a)$.

\begin{definition}
The {\it $s$-length} of the arc $\a$ is
$s(\a)=\cosh(\ell_{\a}/2)$.
\end{definition}

\begin{remark}\label{rem:s-lengths}
The definition above is due to
Ushijima \cite{ushijima:decomposition}
up to a factor $\sqrt{2}$.
\end{remark}

As a triangulation of a hyperbolic surface produces
a dissection into hyperbolic hexagons with geodesic
edges and right angles, we have the following.

\begin{proposition}[\cite{ushijima:decomposition}]
Let $S$ be a compact hyperbolic Riemann surface with boundary.
Fix a triangulation $\ua=\{\a_1,\dots,\a_{6g-6+3n}\}$
of $S$. The map $s(\ua):\Teich(S)\rar \R_+^{6g-6+3n}$ is
a real-analytic diffeomorphism.
\end{proposition}

Every $S$-marked hyperbolic surface $(\Si,g)$ has a preferred
system of arcs, that is actually a triangulation for most surfaces.
It can be equivalently
obtained using the convex
hull construction \cite{ushijima:decomposition} (following
\cite{penner:decorated}, \cite{epstein-penner:decomposition} and
\cite{kojima:decomposition})
or using the spine \cite{bowditch-epstein:triangulations}.
We follow this second way.
%
%
\begin{subsection}{Spine of a Riemann surface with
boundary}\label{sec:spine1}
Let $\Si$ be a compact hyperbolic Riemann surface with
nonempty boundary, and possibly cusps.

We define the {\it valence} $\nu(u)$
of a point $u\in\Si$ which is not
a cusp as the number of shortest geodesics joining $u$ to $\pa\Si$
that realize the distance $d(u,\pa\Si)$.
Clearly $\nu(u)\geq 1$ (the constant geodesic being allowed).

To define the valence of a cusp $c$, consider a geodesic $\g$
ending in $c$. Fix a small embedded horoball $B$ at $c$
and define the reduced length $\ell^{B}_{\g}$ of $\g$ as the length
of the truncated geodesic $\g\setminus B$.
The {\it shortest geodesics} ending at $c$ are the nonconstant geodesics
joining $c$ with $\pa\Si$, which minimize $\ell^{B}$.
Different choices of $B$ change the reduced
length by a constant term, so that shortest geodesics ending
at a cusp are well-defined.
Thus, we can say that the valence $\nu(c)$ of a cusp $c$ is the number
of shortest geodesics ending at $c$.

\begin{center}
\begin{figurehere}
\psfrag{v1}{$v_1$}
\psfrag{v2}{$v_2$}
\psfrag{v3}{$v_3$}
\psfrag{be1}{$\b_{e_1}$}
\psfrag{be2}{$\b_{e_2}$}
\psfrag{be3}{$\b_{e_3}$}
\psfrag{c}{$c$}
\includegraphics[width=0.7\textwidth]{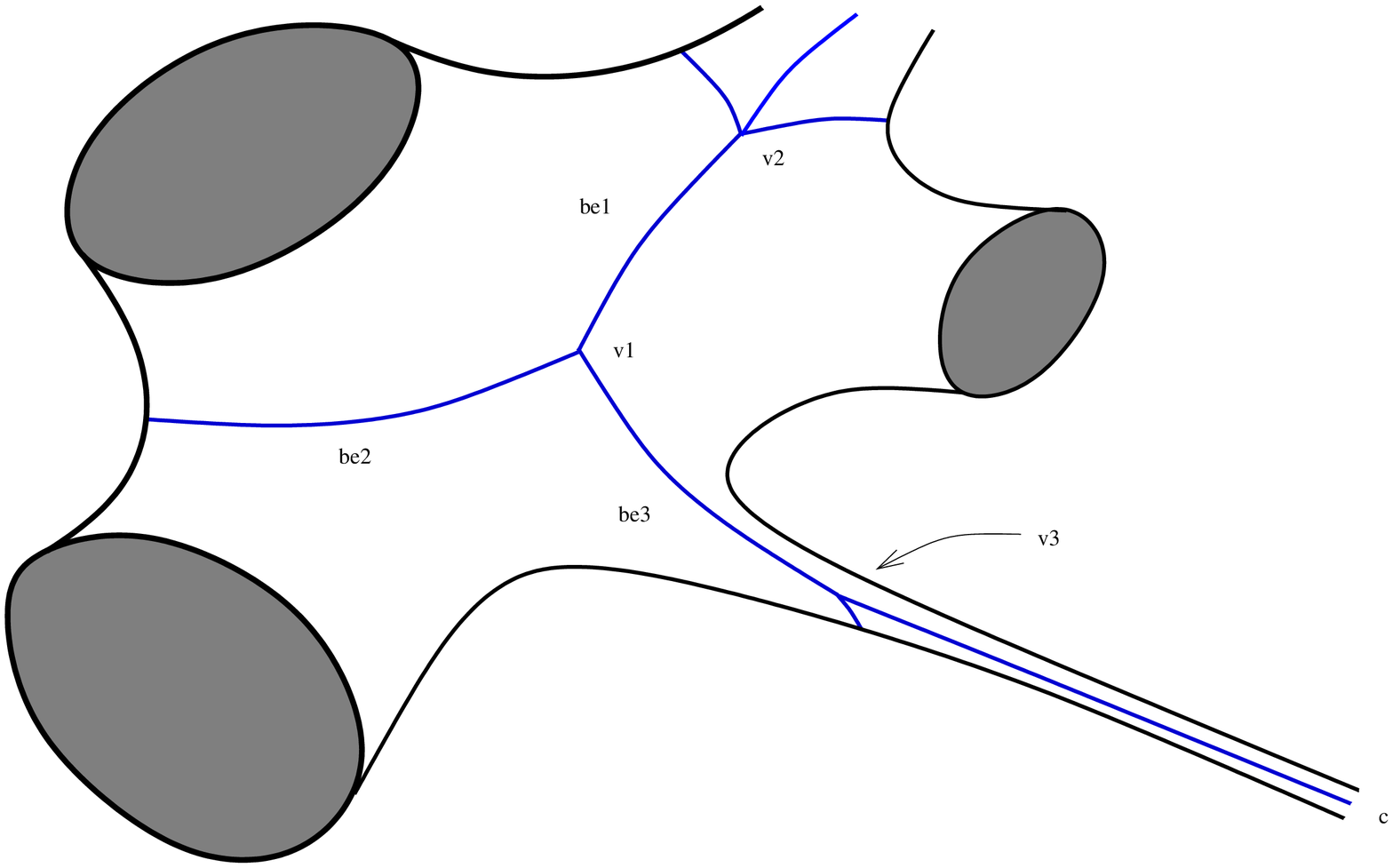}
\end{figurehere}
\end{center}

Define the loci $N:=\{u\in\dot{\Si}\,|\,\nu(u)=2\}$
and $V:=\{u\in\Si\,|\,\nu(u)\geq 3\}\cup\{\text{cusps}\}$.
Notice that $N$ is a finite disjoint
union $N=\coprod_{e\in E_{fin}}\beta_e$
of simple open geodesic arcs ({\it edges})
and $V$ is a finite collection of points ({\it vertices}).

\begin{definition}
The {\it spine} $\Sp(\Si)$
of $\Si$ is the $1$-dimensional CW-complex embedded in
$\Si$ given by $V\cup N$.
\end{definition}

Let $e\in E_{fin}$ be an edge of the spine of $\Si$.
Pick any point $u\in \beta_e$ and let $\g_1$ and $\g_2$
be the two shortest geodesics that join $u$ with $\pa \Si$.
The isotopy class of the unoriented arc with support $\g_1\cup\g_2$
is called {\it dual} to $\b_e$.
We will denote by $\a_e$ the geodesic arc dual to $\b_e$ that meets
the boundary perpendicularly.

Pick $B\subset\Si$ a small embedded horoball at the cusp $c$
such that $B\cap V=\{c\}$ and
call {\it sectors} of the cusp $c$ the connected components of
$B\setminus\Sp(\Si)$.
Clearly, sectors of $c$ bijectively correspond to shortest geodesics
ending in $c$.
Let $E_\infty$ be the set of sectors of all cusps in $\Si$.
For every $e\in E_\infty$, call $\a_e$ the corresponding shortest
geodesic.

Thus, we can attach to $\Si$ a preferred system of arcs
$\Sp^*(\Si):=\Sp^*_{fin}(\Si)\cup\Sp^*_\infty(\Si)\in\Af(\Si)$,
where $\Sp^*_{fin}:=\{\a_e\,|\,e\in E_{fin}\}$
and $\Sp^*_{\infty}:=\{\a_e\,|\,e\in E_\infty\}$.
Call $E:=E_{fin}\cup E_\infty$ and
notice that $\Si\setminus\bigcup_{e\in E}\a_e$ is a disjoint union
of discs, so that $\Sp^*(\Si)$ really belongs to $\Ao(\Si)$.
Also, $\Si\setminus\bigcup_{e\in E_{fin}}\a_e$ is a disjoint union
of discs and pointed discs, so that $\Sp^*_{fin}(\Si)$ belongs
to $\Ao(\Si)$ too.

\begin{center}
\begin{figurehere}
\psfrag{ae'}{$\a_{e'}$}
\psfrag{ae1}{$\a_{e_1}$}
\psfrag{ae2}{$\a_{e_2}$}
\psfrag{ae3}{$\a_{e_3}$}
\psfrag{c}{$c$}
\includegraphics[width=0.8\textwidth]{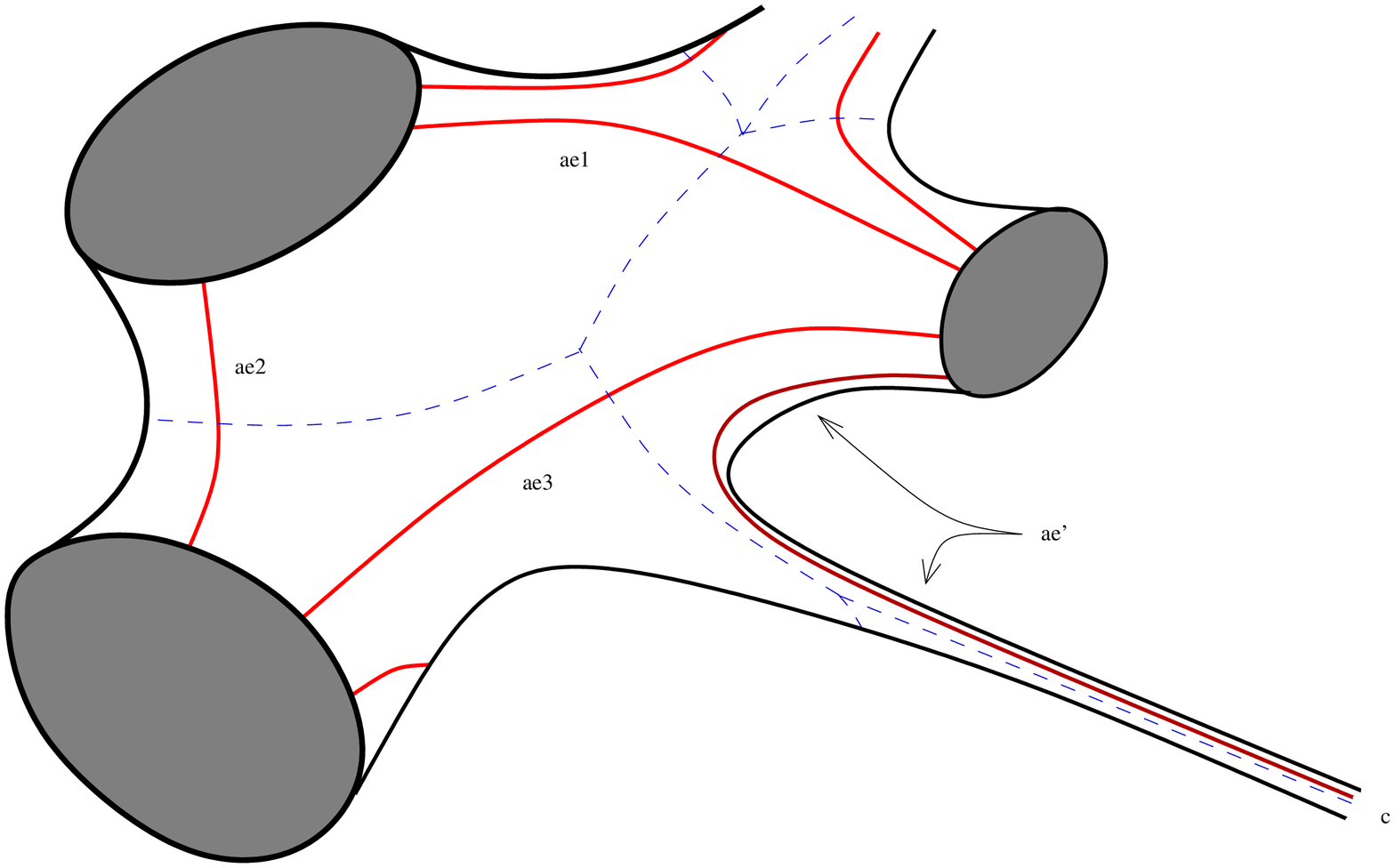}
\caption{$\a_{e'}$ is the shortest geodesic outgoing from
the unique sector of the cusp $c$.}
\end{figurehere}
\end{center}
\end{subsection}
%
%
\begin{subsection}{Spine of a truncated surface}
Let $\Si$ be a compact hyperbolic surface
with boundary circles
$C_1,\dots C_n$ and cusps $c_1,\dots,c_m$.
Fix $\up=(p_1,\dots,p_m)$ a vector of nonnegative real numbers.

For small $\e>0$
let $\Si_{\e\up}$ be the {\it truncated surface} obtained
from $\Si$ by removing the open horoball of radius
$\e p_i$ at the $i$-th cusp (which will be disjoint
for $\e$ small enough).

On $\Si_{\e\up}$ there are well-defined boundary distance
and valence function.
Define $\Sp(\Si_{\e\up})$ to be the spine of $\Si_{\e\up}$.

\begin{lemma}
The closure of
$\lim_{\e\rar 0}\Sp(\Si_{\e\up})$ is a $1$-dimensional
CW-complex embedded in $\Si$, which actually coincides with
the spine $\Sp(\Si)$.
\end{lemma}

\begin{proof}
In fact, for every point $u\in\Si$ which is not a cusp
the restriction of the boundary distance
function $d(-,\pa\Si_{\e\up}):\Si_{\e\up}\rar\R$ to
a fixed small neighbourhood of $u$ stabilizes as $\e\rar 0$
and coincides with the restriction of $d(-,\pa\Si):\Si\rar\R\cup\{\infty\}$.
So that the valence function $\nu_{\e\up}$ also stabilizes.
As $\Sp(\Si_{\e\up})$ is a $1$-dimensional CW-complex for positive $\e$,
the same holds for the limit.
Thus, in this case, the limit is independent of the choice
of $\up$.
\end{proof}

We attach $\Si_{\e\up}$ the system of arcs dual to the spine
$\Sp^*(\Si_{\e\up})\in \Ao(\Si_{\e\up})$.
For $\e$ small, $\Sp^*(\Si_{\e\up})$ coincides with
$\Sp^*(\Si)$, so they define the same arc system in
$\Ao(\Si_{\e\up})\cong\Ao(\Si)$.
\end{subsection}
%
%
\begin{subsection}{Spine of a decorated surface}
Let $\Si$ be a hyperbolic surface with
cusps $c_1,\dots,c_m$ and no boundary circles.
Choose a nonzero
vector of nonnegative numbers $\up=(p_1,\dots,p_m)$ and
denote by $\Si_{\e\up}$ the truncated surface.

As the geodesics that realize the minimum distance from the boundary
meet the horocycles perpendicularly, it is easy to see that
$d(u,\pa\Si_{\e\up})=d(u,\pa\Si_{\e'\up})+\log(\e'/\e)$
for every $u\in\Si_{\e\up}\cap
\Si_{\e'\up}$. Thus, the valence
$\nu$ does not depend on $\e$ (when it is defined), which
essentially proves the following lemma.

\begin{lemma}
The homeomorphism type of $\Sp(\Si_{\e\up})$ stabilizes when $\e\rar 0$.
\end{lemma}

We call $\Sp(\Si,\up)$ the closure
inside $\Si$ of $\lim_{\e\rar 0}\Sp(\Si_{\e\up})$.
Following Section~\ref{sec:spine1}, let $E_{fin}$ be the set of
edges of $\Sp(\Si,\up)$ and $E_\infty$ the set of sectors
of cusps $c_i$ with $p_i=0$. 
Define analogously $\Sp^*(\Si,\up):=\Sp^*_{fin}(\Si,\up)
\cup \Sp^*_\infty(\Si,\up)\in\Ao(\Si)$.

\begin{definition}[\cite{penner:decorated}]\label{def:decorated}
A {\it decorated surface}
is a couple $(\Si,\up)$ where
$\Si$ is a hyperbolic surface with $m$ cusps and no boundary circles,
and $\up=(p_1,\dots,p_m)$ is a nonzero vector of
nonnegative numbers.
\end{definition}

Notice that $\Sp(\Si,\up)$, $\Sp^*(\Si,\up)$
and $\Sp^*_{fin}(\Si,\up)$ depend only on the choice of
a projective class $[\up]\in (\R_{\geq 0}^m\setminus\{0\})/\R_+
\cong \Delta^{m-1}$.
Moreover, the $i$-th cusp is a vertex of the spine
if and only if $p_i=0$.
\end{subsection}
%
%
\begin{subsection}{$\Gamma$-equivariant cellular
decomposition of $\Teich(S)$}\label{sec:decomposition}
Let $\ua$ be a triangulation of a compact hyperbolic
surface $\Si$ with nonempty boundary (and possibly cusps).

Let $\ora{\a_i}$ and
$\ora{\a_j}$ be two distinct oriented
arcs whose supports belong to $\ua$
and which point toward the boundary
component $C$.
Define $d(\ora{\a_i},\ora{\a_j})$ to be the length
of the path along $C$ that runs from the endpoint
of $\ora{\a_i}$ to the endpoint of $\ora{\a_j}$
in the positive direction (according to the orientation
induced on $C$). Clearly, $d(\ora{\a_i},\ora{\a_j})+
d(\ora{\a_j},\ora{\a_i})=\ell_{C}$, which is actually zero
if $C$ is a cusp.

Now, let $\ora{\a_i},\ora{\a_j},\ora{\a_k}$ the oriented
arcs that bound a chosen
connected component $t$ of $\Si\setminus\bigcup_{\a\in\ua}\a$.
Assume that $(\ora{\a_i},\ora{\a_j},\ora{\a_k})$ are cyclically
ordered according to the orientation induced by $t$.
Define
\begin{align*}
& w_{\ua}(\ora{\a_i}):=\frac{1}{2}\left(
d(\ora{\a_i},\ola{\a_j})+d(\ora{\a_k},\ola{\a_i})
-d(\ora{\a_j},\ola{\a_k})\right)
\quad\text{and} \\
& w_{\ua}(\a_i):=w_{\ua}(\ora{\a_i})+w_{\ua}(\ola{\a_i})
\end{align*}
where $\ola{\a}$ is obtained by reversing the orientation of $\ora{\a}$.
\begin{definition}
Given an arc $\a$ of the triangulation $\ua\in\Af(\Si)$,
we call $w_{\ua}(\a)$ the {\it width} of $\a$ with respect to $\ua$.
\end{definition}
\begin{remark}
Luo \cite{luo:decomposition} used the term ``$E$-invariant'' for the width.
\end{remark}
\begin{proposition}\label{prop:width}
With the notation above,
\[
\sinh(w_{\ua}(\ora{\a_i}))=\frac{s(\a_j)^2+s(\a_k)^2-s(\a_i)^2}
{2s(\a_j) s(\a_k)\sqrt{s(\a_i)^2-1}}
\]
\end{proposition}
\begin{proof}
We prove the statement in the case when $w_{\ua}(\ora{\a_r})\geq 0$
for $r=i,j,k$.
The other cases can be treated similarly.
We will denote by $a_r$ the length of $\a_r$
and by $f_r$ the orthogonal projection of $u$
on the side of $H$ facing $\a_r$ for $r=i,j,k$
(see Figure~\ref{fig:hexagon}).
\begin{center}
\begin{figurehere}
\psfrag{zi}{$z_i$}
\psfrag{yi}{$y_i$}
\psfrag{mi}{$m_i$}
\psfrag{fi}{$f_i$}
\psfrag{zj}{$z_j$}
\psfrag{yj}{$y_j$}
\psfrag{mj}{$m_j$}
\psfrag{fj}{$f_j$}
\psfrag{zk}{$z_k$}
\psfrag{yk}{$y_k$}
\psfrag{mk}{$m_k$}
\psfrag{fk}{$f_k$}
\psfrag{u}{$u$}
\psfrag{gi}{$\g_i$}
\psfrag{gj}{$\g_j$}
\psfrag{gk}{$\g_k$}
\includegraphics[width=0.9\textwidth]{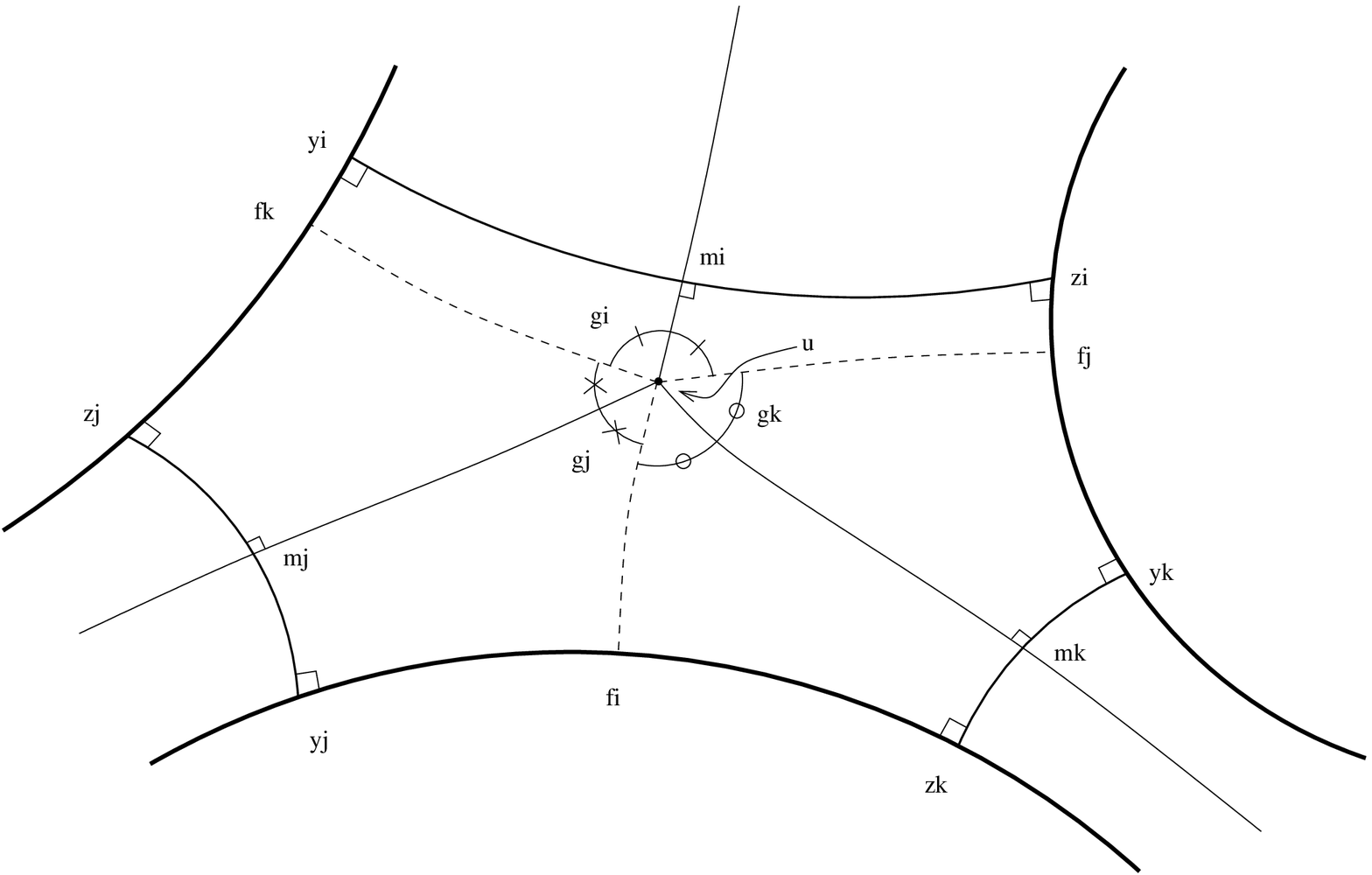}
\caption{Geometry of the hexagon $t$.}
\label{fig:hexagon}
\end{figurehere}
\end{center}
Call $\g_i$ the angle $\wh{m_i\,u\,f_k}=\wh{f_j\,u\,m_i}$
and define analogously $\g_j$ and $\g_k$.
Notice that $\g_i+\g_j+\g_k=\pi$.
As $(m_j\,y_j\,z_k\,m_k\,u)$ is a pentagon with four right angles,
Lemma~\ref{lemma:pentagon} gives
\[
\cosh(\os{\frown}{y_j z_k})=
\frac{\cosh(a_j/2)\cosh(a_k/2)+\cos(\g_j+\g_k)}
{\sinh(a_j/2)\sinh(a_k/2)}
\]
As $(z_i\,y_i\,z_j\,y_j\,z_k\,y_k)$ is an hexagon with six right angles,
by Lemma~\ref{lemma:hexagon} we have
\begin{align*}
\cosh(\os{\frown}{y_j z_k}) & =
\frac{\cosh(a_j)\cosh(a_k)+\cosh(a_i)}
{\sinh(a_j)\sinh(a_k)}=\\
& =\frac{\cosh(a_j)\cosh(a_k)+\cosh(a_i)}
{4\sinh(a_j/2)\sinh(a_k/2)\cosh(a_j/2)\cosh(a_k/2)}
\end{align*}
so that
$\dis\quad
\cosh(a_j/2)\cosh(a_k/2)+\cos(\g_j+\g_k)=
\frac{\cosh(a_j)\cosh(a_k)+\cosh(a_i)}
{4\cosh(a_j/2)\cosh(a_k/2)} .$

As $(m_i\,y_i\,f_k\,u)$ is a quadrilateral with three right angles,
then part (a) of Lemma~\ref{lemma:quadrilateral} gives
\[
\sinh(a_i/2)\sinh(\os{\frown}{y_i f_k})
=\cos(\g_i)=-\cos(\g_j+\g_k)
\]
Because $w_{\ua}(\ora{\a_i})=\os{\frown}{y_i f_k}$,
we deduce
\[
\cosh(a_j/2)\cosh(a_k/2)-\sinh(a_i/2)\sinh(w_{\ua}(\ora{\a_i}))=
\frac{\cosh(a_j)\cosh(a_k)+\cosh(a_i)}
{4\cosh(a_j/2)\cosh(a_k/2)}.
\]

From $s(\a_r)=\cosh(a_r/2)$, we get
$\cosh(a_r)=2s(\a_r)^2-1$ and $\dis\sinh(a_r/2)=\sqrt{s(\a_r)^2-1}$.
Substituting inside the expression above, we get the wanted result.
\end{proof}
\begin{remark}\label{rmk:angle}
As a byproduct of the proof above, we have also obtained that
\[
\cos(\g_i)=\frac{s(\a_j)^2+s(\a_k)^2-s(\a_i)^2}
{2s(\a_j) s(\a_k)}
\]
\end{remark}
Other length functions that are sometimes useful are the
{\it $b$-lengths}: for every hexagon in $\Si\setminus\bigcup_i \a_i$,
the $b$-lengths are the lengths of the edges lying on a boundary
component. In Figure~\ref{fig:hexagon},
the $b$-length $b_{t,i}$ is the length of
the path from $y_j$ to $z_k$ passing through $f_i$.
Using Lemma~\ref{lemma:hexagon}, we have
\[
\cosh(b_{t,i})=\frac{\cosh(a_j)\cosh(a_k)+\cosh(a_i)}{\sinh(a_j)\sinh(a_k)}
\]
Fixed a triangulation, the set of all $b$-lengths is too large
to be a system of coordinates, but Ushijima proved
\cite{ushijima:decomposition} that their relations
are generated by homogeneous quadratic ones in
their hyperbolic cosines.\\

If we deal with the system of arcs $\Sp^*_{fin}(\Si)$ instead of a general
triangulation, we can define the widths even if $\Sp^*(\Si)$ is
not a maximal system.

Consider an arc $\a_e$ with $e\in E_{fin}$,
choose an orientation $\ora{\a_e}$
of $\a_e$ and view $\ora{\a_e}$ as pointing upwards.
For every point $u\in\b_e$, call $P_{\ora{\a_e}}(u)$
the {\it projection} of $u$ to the boundary component
pointed by $\ora{\a_e}$.

Call $w_{sp}(\ora{\a_e})$
the length (with sign)
of the boundary arc that
runs from the endpoint of $\ora{\a_e}$ leftward 
to the projection $P_{\ora{\a_e}}(v_l)$
of the left endpoint $v_l$ of $\b_e$.
Define $w_{sp}(\a_e):=w_{sp}(\ora{\a_e})+w_{sp}(\ola{\a_e})$.

\begin{remark}
The width $w_{sp}(\a_e)$ is always positive, but $w_{sp}(\ora{\a_e})$ or
$w_{sp}(\ola{\a_e})$ might be zero or negative.
Notice that, given $[f:S\rar(\Si,g)]\in\eTeich(S)\setminus\eTeich(S)(0)$
the system $f^* w_{sp}$ defines a point in
$|\Ao(S)|_\R$.
\end{remark}

The following results by Ushijima and Luo adapt and generalize
Penner's work on ideal triangulations
\cite{penner:decorated} (see Section~\ref{sec:decorated}).

\begin{theorem}[\cite{ushijima:decomposition}]
Let $\Si$ be a hyperbolic surface with nonempty boundary
and possibly cusps.
There is at least one triangulation $\ua$ such that
$w_{\ua}(\a)\geq 0$ for all $\a\in\ua$.
Moreover, the intersection of all these triangulations
is $\Sp^*(\Si)$ and $w_{sp}(\a)>0$ for all $\a\in\Sp^*_{fin}(\Si)$.
\end{theorem}

\begin{theorem}[\cite{luo:decomposition}]\label{thm:decomposition}
Let $S$ be a compact hyperbolic surface with boundary.
The induced map
\[
\xymatrix@R=0pt{
\eTeich(S)\setminus\eTeich(S)(0) \ar[r]              &  |\Ao(S)|_\R \\
[f:S\rar(\Si,g)] \ar@{|->}[r] &  f^*w_{sp}
}
\]
is a $\Gamma(S)$-equivariant homeomorphism.
\end{theorem}

Theorem~\ref{thm:decomposition} gives a $\Gamma(S)$-equivariant
cellular decomposition of $\eTeich(S)\setminus\eTeich(S)(0)$
and so an orbisimplicial
decomposition of the moduli space $\eM(S)\setminus\eM(S)(0)$.
\end{subsection}
%
%
\begin{subsection}{The cellular decomposition
for decorated surfaces}\label{sec:decorated}
Let $(\Si,\up)$ be a decorated hyperbolic surface
(see Definition~\ref{def:decorated}) with cusps $c_1,\dots,c_m$
and no boundary, and let $\ua\in\Af(\Si)$
be a triangulation.

Take a small $\e>0$ such that
the horoballs at $c_1,\dots,c_m$
with radii $\e p_1,\dots,\e p_m$ are embedded and disjoint.
The truncated length $\ell^{\e\up}_{\a}$ of an arc
$\a\in\ua$ is
the length of the truncation $\a\cap\Si_{\e\up}$.
As $\ell^{\e'\up}_\a=\ell^{\e\up}_\a+2\log(\e/\e')$
for small $\e,\e'>0$, then we can define
$\ell^{\up}_\a:=\ell^{\e\up}_\a+2\log(\e)$, which is
independent of $\e$.

\begin{theorem}[\cite{penner:decorated}]
Let $S$ be a hyperbolic surface with $m$ cusps
and let
$\ua\in\Af(S)$ be a triangulation.
The lengths $\{\ell^{\up}_{\a}\,|\,\a\in\ua\}$ are
real-analytic coordinates on the space
$\eTeich(S)(0)\times \R_+^m$ of positively decorated
surfaces $([f:S\rar\Si],\up)$.
\end{theorem}

In the theorem above, really Penner used the $\l$-lengths,
defined as $\dis\l(\a,\up):=\sqrt{2\exp(\ell^{\up}_\a)}$.
Notice that Penner's $\l$-lengths are the limit
of the $s$-lengths in the following sense:
given a sequence $([f_n:S\rar\Si_n])$ in $\Teich(S)$
that converges to $[f:S\rar\Si]\in\eTeich(S)(0)$, we have
\[
\lim_{n\rar \infty}
\frac{s(\a_i)(f_n)}{s(\a_j)(f_n)}=
\frac{\l(\a_i,\up^{(\infty)})(f)}{\l(\a_j,\up^{(\infty)})(f)}
\]
whenever $[\up^{(n)}]\rar[\up^{(\infty)}]$ in $\Delta^{m-1}$
(and we have set $\up^{(n)}=\Ll(f_n)$).

On the other hand, the role
of the distance $d(\ora{\a_i},\ora{\a_j})$
between the endpoints of two oriented
arcs $\ora{\a_i}$ and $\ora{\a_j}$ (defined in
Section~\ref{sec:decomposition}) is played by
the length $d_{\up}(\ora{\a_i},\ora{\a_j})$ of the
horocyclic segment,
running from the endpoint of $\ora{\a_i}\cap\Si_{\up}$
to the endpoint of $\ora{\a_j}\cap\Si_{\up}$ in the positive direction.

Pick a connected component $t$ of $\Si\setminus\bigcup_{\a\in\ua}\a$
and let $(\ora{\a_i},\ora{\a_j},\ora{\a_k})$ the arcs that bound $t$
with the induced orientation and cyclic order.
%

%
Penner defined the ``simplicial coordinate''
$X_i:=X(\ora{\a_i},\up)+X(\ola{\a_i},\up)$ associated to $\a_i$
setting
\[
X(\ora{\a_i},\up):=\frac{\l(\a_j,\up)^2+\l(\a_k,\up)^2-\l(\a_i,\up)^2}
{\l(\a_i,\up)\l(\a_j,\up)\l(\a_k,\up)}
\]
For a sequence $([f_n:S\rar\Si_n])$ as above,
we also have
\[
\lim_{n\rar\infty} \frac{w_{\ua}(\ora{\a_i})(f_n)}
{\sum_{k=1}^m p_k^{(n)}}=\frac{X(\ora{\a_i},\up^{(\infty)})(f)}
{\sum_{k=1}^m p_k^{(\infty)}}
\]
Similarly, fixed a triangulation $\ua$,
Penner defined the {\it $h$-lengths} to be the lengths
one-half the lengths of the horocyclic arcs appearing in
the truncated triangles of $\Si_{\e\up}\setminus\bigcup_i\a_i$.
If $t$ is a truncated triangle, bounded by arcs $\a_i$, $\a_j$,
$\a_k$ (cyclically ordered), then Penner showed that
$\dis h_{t,i}=\frac{\lambda_i}{\lambda_j \lambda_k}$.
One can observe that
\[
\lim_{n\rar\infty}\frac{b_{t,i}(f_n)}{\sum_{k=1}^m p_k^{(n)}}=
\frac{2h_{t,i}(f,\up^{(\infty)})}{\sum_{k=1}^m p_k^{(\infty)}}
\]
so the $b$-lengths limit to the $h$-lengths (up to a factor $2$).

The convex hull construction, or equivalently the spine
$\Sp(\Si,\up)$,
gives a preferred system of arcs $\Sp^*(\Si,\up)$ on $\Si$.
Analogously to what done in Section~\ref{sec:decomposition}
with the widths,
one can define simplicial coordinates $X_{sp}$
for arcs in $\Sp^*(\Si,\up)$ as half the
lengths of their projection to the
truncating horocycles.

\begin{theorem}[\cite{penner:decorated}]
Let $(\Si,\up)$ be a hyperbolic decorated surface.
There is at least one triangulation $\ua$ such that
$X(\a,\up)\geq 0$ for all $\a\in\ua$.
Moreover, the intersection of all these triangulations
is $\Sp^*(\Si)$ and $X_{sp}(\a,\up)>0$ for all $\a\in\Sp^*_{fin}(\Si)$.
\end{theorem}

\begin{theorem}[\cite{penner:decorated}]\label{thm:decomposition-decorated}
Let $S$ be a hyperbolic surface with $m$ boundary components.
The induced map
\[
\xymatrix@R=0pt{
\eTeich(S)(0)\times (\Delta^{m-1}\times\R_+) \ar[r] &  |\Ao(S)|_\R \\
([f:S\rar(\Si,g)],\up) \ar@{|->}[r] &  f^* X_{sp}
}
\]
is a $\Gamma(S)$-equivariant homeomorphism.
\end{theorem}

Theorem~\ref{thm:decomposition-decorated} provides a $\Gamma(S)$-equivariant
cellular decomposition of $\eTeich(S)(0)\times(\Delta^{m-1}\times\R_+)$.
and so an orbisimplicial
decomposition of $\eM(S)(0)\times (\Delta^{m-1}\times\R_+)$.
\end{subsection}
%
%
\end{section}
%
%
\begin{section}{Fenchel-Nielsen deformations and distances between
geodesics}\label{sec:FN}
%
%
\begin{subsection}{The Fenchel-Nielsen deformation}\label{sec:FN-deformation}
Let $R$ be a hyperbolic surface without
boundary and $\xi\subset R$ a simple closed geodesic.
A (right) {\it Fenchel-Nielsen deformation} $\mathrm{Tw}_{t\xi}$
of $R$ along $\xi$ of translation distance $t$
is obtained by cutting $R$ along $\xi$,
sliding the left side forward by $t$
relatively to the right side and regluing the two sides.
Notice that the deformation is an isometry outside $\xi$.

The terminology is due to the fact that the deformation
pushes one to the right when one passes the default line.

\begin{remark}
Let $R$ be a compact hyperbolic surface without boundary
and $\{\xi_1,\dots,\xi_N\}$
a maximal system of simple closed curves. Let $(\ell_i,\tau_i)_{i=1}^N$ be
the associated Fenchel-Nielsen coordinates.
Then, the Fenchel-Nielsen deformation along $\xi_1$ of translation
distance $t$ acts as
$(\ell_1,\tau_1,\dots,\ell_N,\tau_N)\mapsto
(\ell_1,\tau_1+t,\ell_2,\tau_2,\dots,\ell_N,\tau_N)$.
\end{remark}

The Fenchel-Nielsen deformation
$\Tw_{t\xi}:\Teich(R)\rar\Teich(R)$
is the flow of a {\it Fenchel-Nielsen vector field}
$\pa/\pa\tau_\xi$ on $\Teich(R)$ (see \cite{wolpert:symplectic}).\\

Let $\H=\{z\in\C\,|\,\mathrm{Im}(z)>0\}$ be the Poincar\'e upper half-plane
and let $\pa\H=\ol{\R}=\R\cup\{\infty\}$ the extended real line.

Choose a uniformization $\pi:\H\rar R$ and let $G=\aut(\pi)
\cong\pi_1(R)$.
Fix a simple closed geodesic $\xi\subset R$ and
let $\ti{\xi}=\os{\frown}{s_1 s_2}\subset\H$
be a {\it lift} of $\xi$, that is a connected component of $\pi^{-1}(\xi)$,
where $s_1,s_2\in\pa\H$ and
$\os{\frown}{s_1 s_2}$ denotes the geodesic on $\H$ with limit
points $s_1$ and $s_2$.
The lift of $\Tw_{t\xi}$ is the composition of the Fenchel-Nielsen
deformations $\Tw_{t\tilde{\xi}}$ along all the lifts
$\ti{\xi}$ of $\xi$.

The Fenchel-Nielsen deformation of $R\cong\H/G$
can be described
as $\H/w_t G w_t^{-1}$, where $(w_t:\H\rar\H)_t$ is a continuous family of
quasi-conformal automorphisms
that fix $0$, $1$ and $\infty$, and $w_0$ is the identity.

A typical case (described in \cite{wolpert:symplectic})
is when $G$ is the cyclic group generated by
the hyperbolic transformation $(z\mapsto \l z)$ with $\l>0$
and the Fenchel-Nielsen deformation is performed along the
simple closed geodesic $\pi(\os{\frown}{0\infty})$.

Let $\theta=\arg(z)$ and $\Phi(\theta)=\int_0^\theta\varphi(\alpha)d\alpha$,
where $\varphi:(0,\pi)\rar\R_{\geq 0}$ is a smooth function with
compact support and $\int_0^\pi \varphi(\alpha) d\alpha=1/2$.
Then, $w_t$ is given by
\begin{equation}\label{eq:wt}
w_t(z)=z\cdot\exp[2 t\Phi(\theta)]
\end{equation}

In this case,
by $\pa z/\pa\tau_{\os{\frown}{0\infty}}$ (at the identity) we will mean
$\pa w_t(z)/\pa t$ (evaluated at $t=0$) for every $z\in\ol{\H}$.
\end{subsection}
%
%
\begin{subsection}{Cross-ratio and Fenchel-Nielsen deformation}
Endow the extended real line $\ol{\R}=\R\cup\{\infty\}$
with the natural cyclic ordering $\prec$ coming from
$\ol{\R}\cong S^1$.
Given $p,q,r,s\in\ol{\R}$,
their {\it cross-ratio} $(p,q,r,s)\in\ol{\R}$ is defined as
\[
(p,q,r,s):=\frac{(p-r)(q-s)}{(p-s)(q-r)}
\]
Wolpert computed how the cross-ratio $(p,q,r,s)$ varies
under infinitesimal Fenchel-Nielsen deformation of $\H$ along
the geodesic $\os{\frown}{s_1 s_2}$ with limit points $s_1,s_2\in\ol{\R}$.

\begin{lemma}[\cite{wolpert:symplectic}]\label{lemma:cr}
Assume $z_1,z_2,z_3,z_4\in\ol{\R}$ are distinct
and $s_1,s_2\in\ol{\R}$ are distinct.
Then
\[
\frac{\pa}{\pa\tau_{\os{\frown}{s_1 s_2}}}(z_1,z_2,z_3,z_4)=
(z_1,z_2,z_3,z_4)\sum_{j=1}^4 \chi_L(z_j)
\left[
(z_{\s(j)},s_1,s_2,z_j)-(z_{\tau(j)},s_1,s_2,z_j)
\right]
\]
where $\dis\s=\left(\begin{array}{cccc} 1 & 2 & 3 & 4 \\
3 & 4 & 1 & 2\end{array}\right),
\ \tau=\left(\begin{array}{cccc} 1 & 2 & 3 & 4 \\
4 & 3 & 2 & 1\end{array}\right)\in\mathfrak{S}_4$
and $\chi_L$ is the characteristic function of $[s_2,s_1]\subset\ol{\R}$
(where $[s_2,s_1]=\{x\in\ol{\R}\,|\,s_1\prec x\prec s_2\}$).
\end{lemma}
The proof follows from the explicit expression of $w_t$
given in Equation~\ref{eq:wt}.
\begin{multicols}{2}{
Consider two nonintersecting geodesics $\os{\frown}{qr}$ and
$\os{\frown}{ps}$
in the upper half-plane $\H$ with endpoints $p,q,r,s\in\ol{\R}$.
The distance $h=\ell_\delta$ between $\os{\frown}{qr}$ and
$\os{\frown}{ps}$
is given by $\cosh (h)=1-2(p,q,r,s)$, or equivalently
$(p,q,r,s)=-\sinh^2 (h/2)$.

\begin{figurehere}
\psfrag{d}{$\delta$}
\psfrag{p}{$p$}
\psfrag{q}{$q$}
\psfrag{r}{$r$}
\psfrag{s}{$s$}
\psfrag{H}{$\H$}
\includegraphics[width=0.4\textwidth]{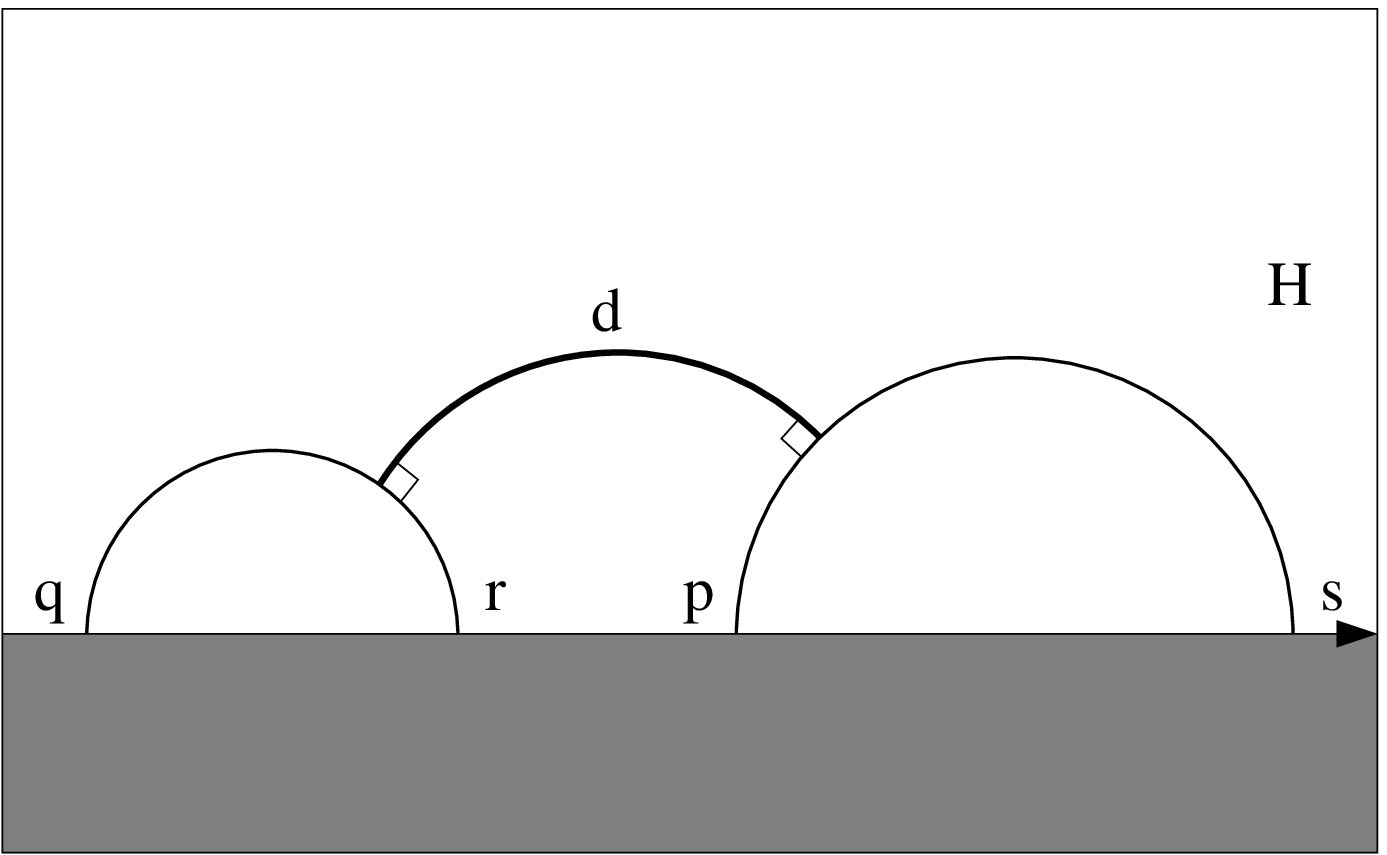}
\end{figurehere}
}
\end{multicols}
\begin{multicols}{2}{
\begin{figurehere}
\psfrag{e}{$\varepsilon$}
\psfrag{x}{$x$}
\psfrag{p}{$p$}
\psfrag{q}{$q$}
\psfrag{r}{$r$}
\psfrag{s}{$s$}
\psfrag{H}{$\H$}
\includegraphics[width=0.4\textwidth]{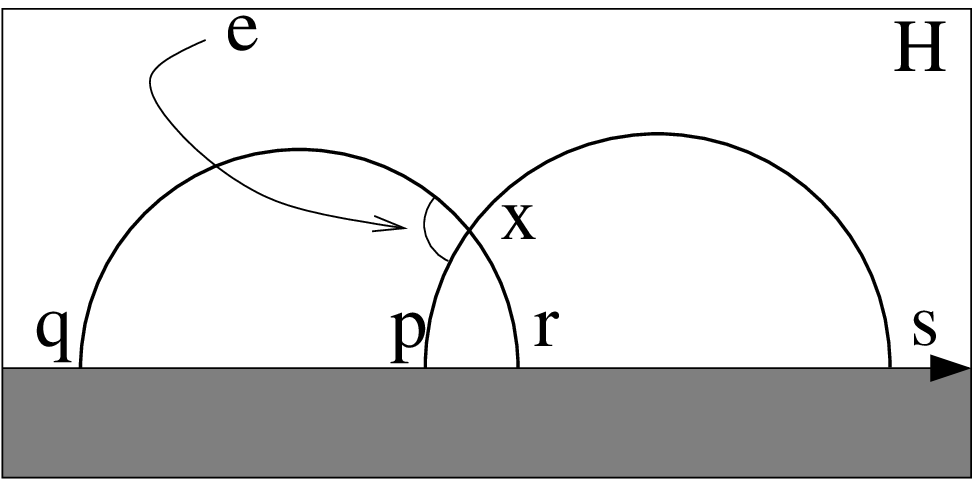}
\end{figurehere}

{
Consider two intersecting geodesics $\os{\frown}{qr}$ and
$\os{\frown}{ps}$
in the upper half-plane $\H$ with endpoints $p,q,r,s\in\ol{\R}$.
The angle $\e=\wh{pxq}$
is given by $\cos (\e)=2(p,q,r,s)-1$, or equivalently
$(p,q,r,s)=\cos^2 (\e/2)$.}
}
\end{multicols}
\end{subsection}
%
%
\begin{subsection}{The variation of the distance between two geodesics}
Let $R$ be a hyperbolic surface without boundary and let
$\g_1$ and $\g_2$ are two (possibly closed) geodesics in $R$.
Let $\delta\subset R$ a (nonconstant) geodesic arc meeting $\g_1$ and $\g_2$
perpendicularly at its endpoints $y_1,y_2$.
Orient $\g_i$ in such a way that, if we travel along $\g_i$ in 
the positive direction, then at $y_i$ we see $\delta$ on our left.

Let $\xi\subset R$ be a simple closed geodesic.
If $x_i\in \xi\cap\g_i$, then we will denote by $\nu_{x_i}$ the
positive angle
at $x_i$ formed by a positively oriented vector along $\g_i$ and $\xi$
and by $d(y_i,x_i)$ the
length
of the path obtained by travelling from $y_i$ to $x_i$
along $\g_i$ (which is a well-defined real number, if $\g_i$ is open,
whereas it is required to belong
to the interval $(0,p_i)$, if $\g_i$ is closed).
The proof of the theorem below adapts arguments of Wolpert
in \cite{wolpert:symplectic}.

\begin{theorem}\label{thm:derivative}
With the above notation, assume that $\xi$ are $\delta$
are disjoint and that $\nu_x=\pi/2$ for every
$x\in\xi\cap(\g_1\cup \g_2)$.
Then
\[
\frac{\pa}{\pa\tau_{\xi}}(h)=c_1+c_2
\qquad\text{with}\quad c_i=
\sum_{x_i\in\xi\cap\g_i}
c_i(x_i)
\]
where $p_i=\ell_{\g_i}$, $h=\ell_\d$ and
\[
c_i(x_i)=
\begin{cases}
\dis\frac{\s}{2}\exp(-|d(y_i,x_i)|) & \text{if $\g_i$ is open}\\
\phantom{x}\\
\dis\frac{\sinh(p_i/2-d(y_i,x_i))}{2\sinh(p_i/2)}
& \text{if $\g_i$ is closed}
\end{cases}
\]
with $\s=\mathrm{sgn}(d(y_i,x_i))$.
\end{theorem}
\begin{remark}
To check that the result and the coefficients
in the formula above are reasonable, think of the
case when $\xi$ intersects $\g_1$ in one point
$x_1$ and it does not intersect $\g_2$.
If $d(y_1,x_1)>0$ is very small, then
the derivative is close to $1/2$. In fact, after
performing a right twist of length $\e$ along $\xi$,
the new geodesic $\g_1$ will interpolate the two broken
branches of the old $\g_1$, and so it will be farther
from $\g_2$ by $\e/2$.
\end{remark}

Choose $\pi:\H\rar R$ a uniformization
and pick a lift $\ti{\d}\subset\H$ of $\d$.
Call $\ti{y}_i$ the endpoint of $\ti{\d}$ mapped to $y_i$.
Let $\ti{\g}_i$ be the lift of $\g_i$ passing through $\ti{y}_i$
and call $p,q,r,s\in\ol{\R}$ their ideal endpoints in such
a way that $\ti{\g}_1=\os{\frown}{ps}$, $\ti{\g}_2=\os{\frown}{qr}$
and $p\prec s\prec q\prec r\prec p$ in the cyclic order $\prec$ of
$\ol{\R}\cong S^1$.
Call the portions $\ti{\g}_1^+:=\os{\frown}{\ti{y}_1 s}\subset\ti{\g}_1$
(resp. $\ti{\g}_1^-:=\os{\frown}{p\ti{y}_1}\subset\ti{\g}_1$) and
$\ti{\g}_2^+:=\os{\frown}{\ti{y}_2 r}\subset\ti{\g}_2$
(resp. $\ti{\g}_2^-:=\os{\frown}{q\ti{y}_2}\subset\ti{\g}_2$)
{\it positive} (resp. {\it negative}).
Under the hypotheses of Theorem~\ref{thm:derivative},
a lift $\ti{\xi}$ of $\xi$ does not intersect $\ti{\delta}$
and it may intersect at most
one of the four geodesic segments $\ti{\g}_i^{\pm}$.
\begin{center}
\begin{figurehere}
\psfrag{de}{$\sst\tilde{\d}$}
\psfrag{p}{$\sst p$}
\psfrag{q}{$\sst q$}
\psfrag{r}{$\sst r$}
\psfrag{s}{$\sst s$}
\psfrag{s1}{$\sst s_1$}
\psfrag{s2}{$\sst s_2$}
\psfrag{x1}{$\sst \ti{x}_1^k$}
\psfrag{x2}{$\sst \ti{x}_2$}
\psfrag{y1}{$\sst \ti{y}_1$}
\psfrag{y2}{$\sst \ti{y}_2$}
\psfrag{g1-}{$\sst \ti{\g}_1^-$}
\psfrag{g1+}{$\sst \ti{\g}_1^+$}
\psfrag{g2-}{$\sst \ti{\g}_2^-$}
\psfrag{g2+}{$\sst \ti{\g}_2^+$}
\psfrag{n}{$\sst \nu_{x_1}$}
\psfrag{xi}{$\sst\ti{\xi}$}
\psfrag{km0}{$k>0$}
\includegraphics[width=0.5\textwidth]{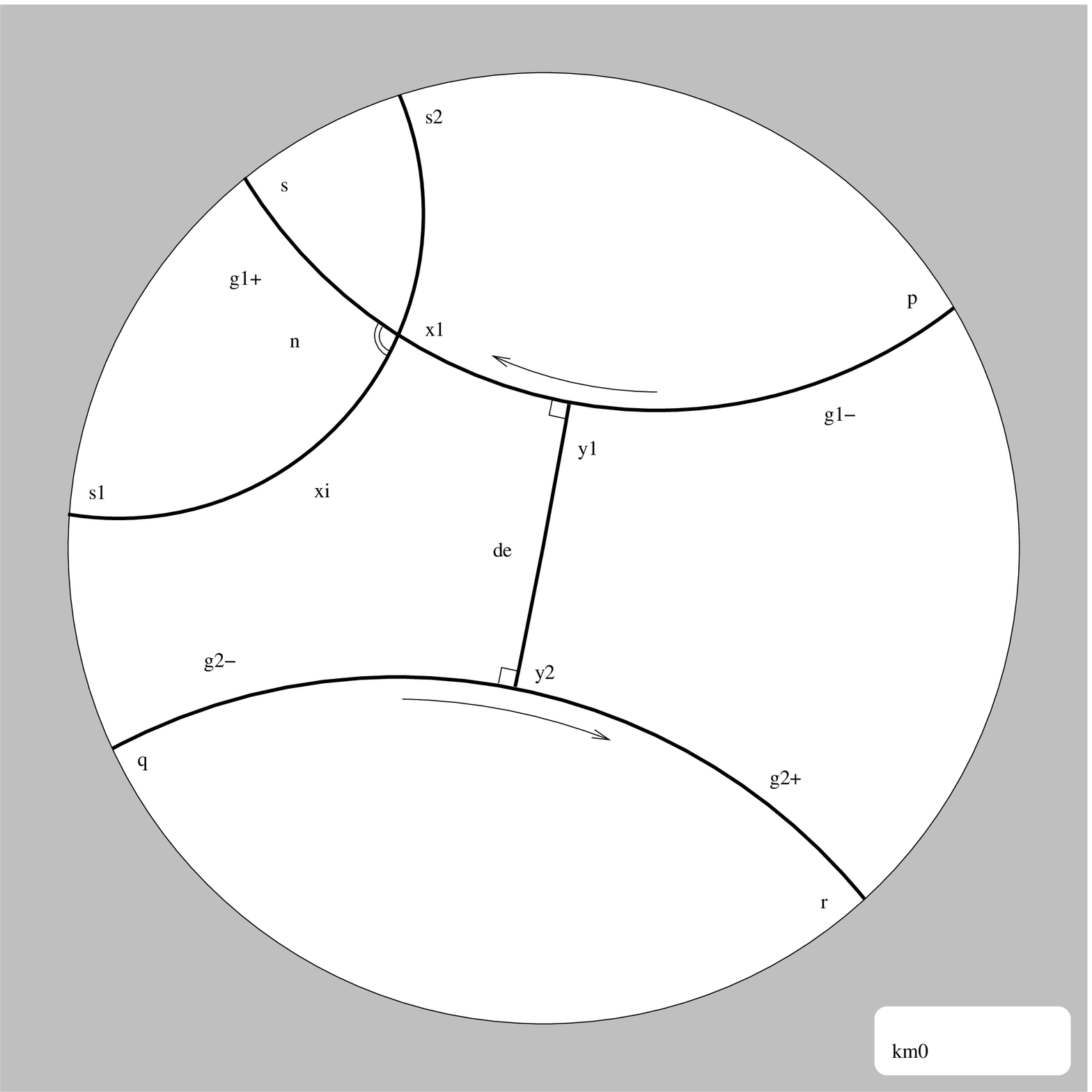}
\end{figurehere}
\end{center}
%
Pick $x_i\in\g_i\cap\xi$. If $\g_i$ self-intersects at $x_i$,
then consider each branch of $\g_i$ separately.

If $\g_i$ is open, then
$x_i$ has only one lift that lies on $\tilde{\g}_i$.
Call it $\ti{x}^0_i$ if it belongs to $\ti{\g}_i^+$
and $\ti{x}^{-1}_i$ if it belongs to $\ti{\g}_i^-$.

If $\g_i$ is closed, then
let $\ti{x}^0_i$ the lift of $x_i$
which belongs to $\ti{\g}_i^+$ and which is closest to $\ti{y}_i$.
Consider $\g_i$ as a loop based at $x_i$ and define
$\ti{x}^k_i$ to be the endpoint
of the lift of $(\g_i)^k$ that starts at $\ti{x}^0_i$ for
every $k\in\mathbb{Z}$.
Clearly, $\ti{x}^k_i\in\ti{\g}_i^+$ for $k\geq 0$ and
$\ti{x}^k_i\in\ti{\g}_i^-$ for $k<0$.
Notice that the distance with sign
$d(\ti{y}_i,\ti{x}^k_i)$ (that is, the length of the portion
of $\ti{\g}_i$ running from $\ti{y}_i$ to $\ti{x}^k_i$
in the positive direction) is exactly $d(y_i,x_i)+kp_i$.

Call $\ti{\xi}(\ti{x}^k_i)$ the only lift of $\xi$ that
passes through $\ti{x}^k_i$.
The derivative $\pa (p,q,r,s)/\pa\tau_\xi$, which we
will sometimes denote by $\pa_\xi(p,q,r,s)$ for brevity, is the sum of
$\pa_{\ti{\xi}}(p,q,r,s)$ for all lifts $\ti{\xi}$ of $\xi$.
Notice immediately that the deformation along $\ti{\xi}$
does not contribute if $\ti{\xi}$ does not intersect
$\ti{\g}_1\cup\ti{\g}_2$.

Define the contribution of $x_i\in\xi\cap\g_i$ to
$\pa(p,q,r,s)/\pa\tau_\xi$ as
\[
\frac{\pa(p,q,r,s)}{\pa\tau_{\xi,x_i}}:=
\begin{cases}
\dis\frac{\pa(p,q,r,s)}{\pa\tau_{\ti{\xi}(\ti{x}_i^0)}}
& \text{if $\g_i$ is open and $d(y_i,x_i)>0$} \\
\dis\frac{\pa(p,q,r,s)}{\pa\tau_{\ti{\xi}(\ti{x}_i^{-1})}}
& \text{if $\g_i$ is open and $d(y_i,x_i)<0$} \\
\dis\sum_{k\in\Z}\frac{\pa(p,q,r,s)}{\pa\tau_{\ti{\xi}(\ti{x}_i^k)}}
& \text{if $\g_i$ is closed}
\end{cases}
\]
and $\dis c_i(x_i)=\frac{\tanh(h/2)}{(p,q,r,s)}\,
\frac{\pa (p,q,r,s)}{\pa\tau_{\xi,x_i}}$.
\end{subsection}
%
%
\begin{subsection}{Proof of Theorem~\ref{thm:derivative}}
Now we compute the contribution of the Fenchel-Nielsen
infinitesimal deformation along $\ti{\xi}(\ti{x}^k_i)$
for $x_i\in\xi\cap\g_i$.
%
%
\begin{subsubsection}{Contribution of $\ti{\xi}(\ti{x}^k_1)$ for $k\geq 0$}
Let $\ti{\xi}(\ti{x}^k_1)=\os{\frown}{s_1 s_2}$ in such a way that
$s\in (s_2,s_1)$ and call $D_k$ the geodesic segment joining $\ti{y}_1$
and $\ti{x}_1^k$.

Lemma~\ref{lemma:cr} gives us
\[
\pa_{\os{\frown}{s_1 s_2}}(p,q,r,s)=(p,q,r,s)[(q,s_1,s_2,s)-(p,s_1,s_2,s)]
\]
The geodesics $\ti{\g}_1=\os{\frown}{ps}$ and
$\os{\frown}{s_1 s_2}$ intersect orthogonally.
Hence, $(p,s_1,s_2,s)=1/2$ and so
$(q,s_1,s_2,s)-(p,s_1,s_2,s)=1/2\cos(\pi-\th)$.
We have so far obtained
\[
\pa_{\os{\frown}{s_1 s_2}}(p,q,r,s)=\frac{1}{2}(p,q,r,s)\cos(\pi-\th)
=-\frac{1}{2}(p,q,r,s)\cos(\th)
\]
where $\th$ is the angle shown in the picture below.

\begin{center}
\begin{figurehere}
\psfrag{de}{$\tilde{\d}$}
\psfrag{p}{$p$}
\psfrag{q}{$q$}
\psfrag{r}{$r$}
\psfrag{s}{$s$}
\psfrag{n}{$n$}
\psfrag{la}{$\l$}
\psfrag{m}{$m$}
\psfrag{th}{$\th$}
\psfrag{s1}{$ s_1$}
\psfrag{s2}{$ s_2$}
\psfrag{x}{$ \ti{x}_1^k$}
\psfrag{y1}{$ \ti{y}_1$}
\psfrag{y2}{$ \ti{y}_2$}
\psfrag{g1}{$ \ti{\g}_1$}
\psfrag{g2}{$ \ti{\g}_2$}
\includegraphics[width=0.6\textwidth]{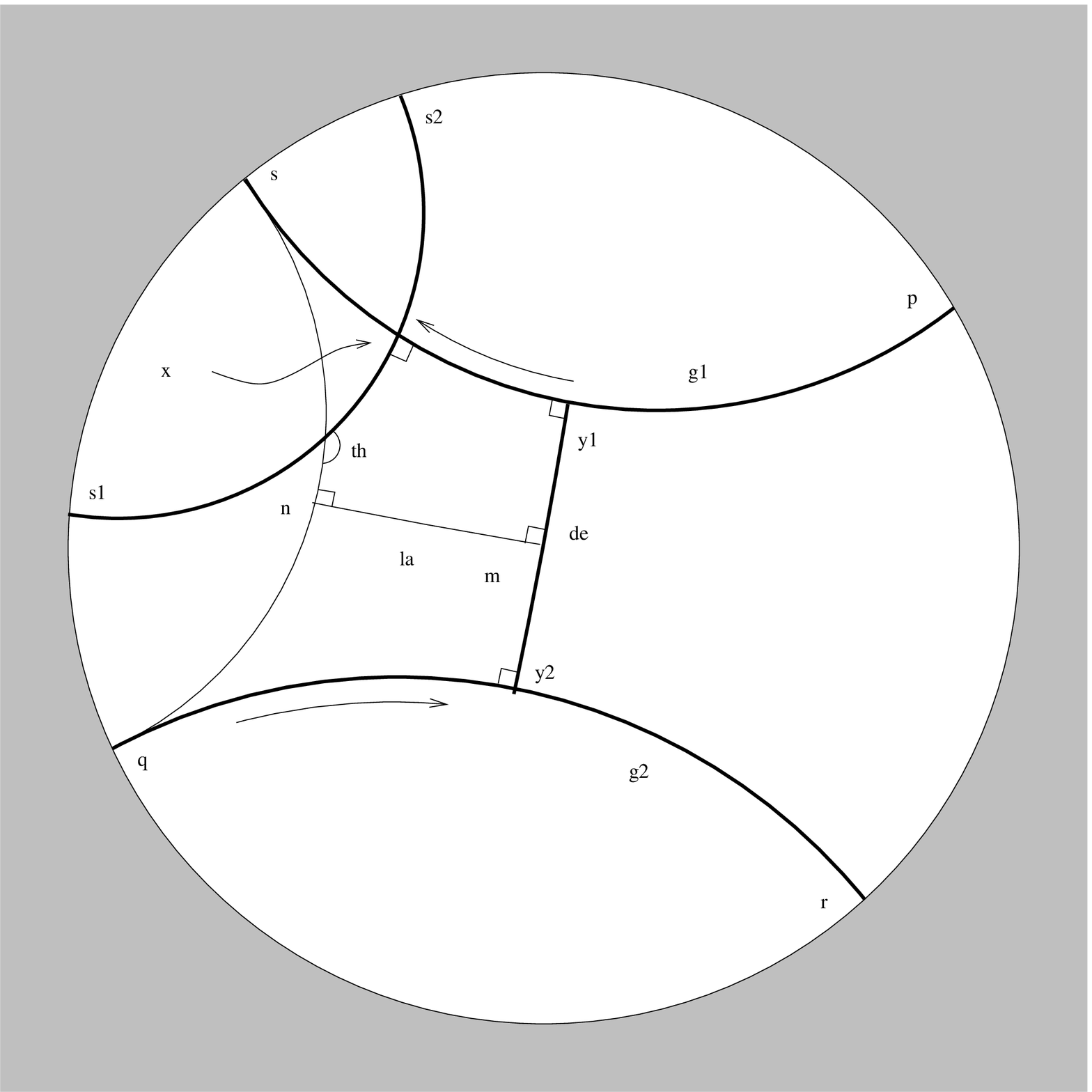}
\caption{Picture for $k\geq 0$}
\end{figurehere}
\end{center}

Let $m$ be the midpoint of $\ti{\delta}$ and 
let $\l$ be the geodesic segment that meets
$\ti{\d}$ and $\os{\frown}{qs}$ orthogonally.
Notice that, if $p,q,r,s$ are fixed,
then $\os{\frown}{s_1 s_2}$ is uniquely determined
by $\ell_{D_k}$ and
$\cos(\th)$ is a real-analytic function of
$\ell_{D_k}$.

In the picture above, we are assuming that
$\os{\frown}{s_1 s_2}$ does not meet $\l$, but
the formula for $\cos(\th)$
we will derive in this case will hold even in
the case when $\os{\frown}{s_1 s_2}$ intersects $\l$,
because of the real-analyticity mentioned above.

Applying part (a) of Lemma~\ref{lemma:quadrilateral}
to the quadrilateral
$(q\,\ti{y}_2\, m\, n)$, we obtain
\[
\sinh(l)=\frac{1}{\sinh(h/2)} \implies
\cosh(l)=\sqrt{1+\frac{1}{\sinh^2(h/2)}}=\frac{1}{\tanh(h/2)}
\]
where $h=\ell_{\ti{\d}}$ and $l=\ell_\l$.

Applying Lemma~\ref{lemma:pentagon}
to the pentagon
$(\th\,n\,m\,\ti{y}_1\,\ti{x}_1^k)$, we obtain
\[
\cosh(h/2)=\frac{\cosh(l)\cosh(d_k)+\cos(\th)}
{\sinh(l)\sinh(d_k)}
\]
where $d_k=\ell_{D_k}$. Thus
\begin{align*}
\cos(\th) = & \cosh(h/2)\sinh(l)\sinh(d_k)-\cosh(l)\cosh(d_k) =\\
= & \frac{\sinh(d_k)-\cosh(d_k)}{\tanh(h/2)}=
 -\frac{\exp(-d_k)}{\tanh(h/2)}
\end{align*}
Hence, $\dis\pa_{\ti{\xi}(\ti{x}_1^k)}(p,q,r,s)
=(p,q,r,s)\frac{\exp(-d_k)}{2\tanh(h/2)}$.
If $\g_1$ is closed, $d_k=d_0+k p_1$ for $k\geq 0$,
and so $\dis\pa_{\ti{\xi}(\ti{x}_1^k)}(p,q,r,s)
=(p,q,r,s)\frac{\exp(-d_0-kp_1)}{2\tanh(h/2)}$.
\end{subsubsection}
%
%
\begin{subsubsection}{Contribution of $\ti{\xi}(\ti{x}^k_1)$ for $k<0$}
Let $\ti{\xi}(\ti{x}^k_1)=\os{\frown}{s_1 s_2}$ in such a way that
$p\in (s_1,s_2)\subset\ol{\R}$.
Lemma~\ref{lemma:cr} gives us
\[
\pa_{\os{\frown}{s_1 s_2}}(p,q,r,s)=(p,q,r,s)[(p,s_1,s_2,r)-(p,s_1,s_2,s)]
\]
As in the previous case, $(p,s_1,s_2,s)=1/2$ and so
$(p,s_1,s_2,r)-(p,s_1,s_2,s)=1/2\cos(\th)$. Thus
\[
\pa_{\os{\frown}{s_1 s_2}}(p,q,r,s)=\frac{1}{2}(p,q,r,s)\cos(\th)
\]
where $\th$ is the angle shown in the picture below.

\begin{center}
\begin{figurehere}
\psfrag{de}{$\tilde{\d}$}
\psfrag{p}{$p$}
\psfrag{q}{$q$}
\psfrag{r}{$r$}
\psfrag{s}{$s$}
\psfrag{n}{$n$}
\psfrag{la}{$\l$}
\psfrag{m}{$m$}
\psfrag{th}{$\th$}
\psfrag{s1}{$ s_1$}
\psfrag{s2}{$ s_2$}
\psfrag{x}{$ \ti{x}_1^k$}
\psfrag{y1}{$ \ti{y}_1$}
\psfrag{y2}{$ \ti{y}_2$}
\psfrag{g1}{$ \ti{\g}_1$}
\psfrag{g2}{$ \ti{\g}_2$}
\includegraphics[width=0.6\textwidth]{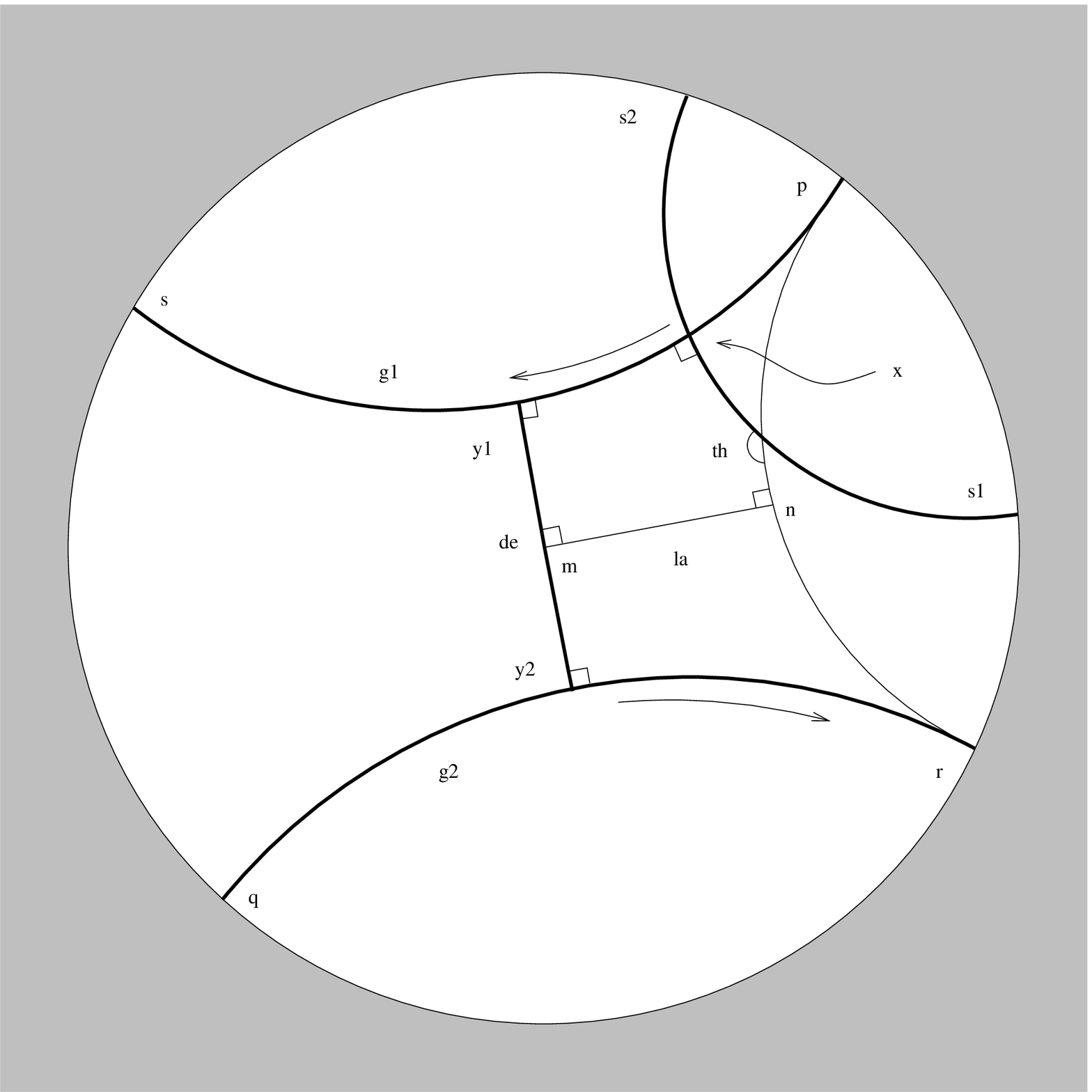}
\caption{Picture for $k<0$}
\end{figurehere}
\end{center}

Arguing as in the case $k\geq 0$, we obtain
\begin{align*}
\cos(\th) = & \cosh(h/2)\sinh(l)\sinh(d_k)-\cosh(l)\cosh(d_k) =\\
= & \frac{\sinh(d_k)-\cosh(d_k)}{\tanh(h/2)}=
 -\frac{\exp(-d_k)}{\tanh(h/2)}
\end{align*}
Hence, $\dis\pa_{\ti{\xi}(\ti{x}_1^k)}(p,q,r,s)
=-(p,q,r,s)\frac{\exp(-d_k)}{2\tanh(h/2)}$.
If $\g_1$ is closed, $d_k=-kp_1-d_0$ for $k<0$ and so
$\dis \pa_{\ti{\xi}(\ti{x}_1^k)}(p,q,r,s)
=-(p,q,r,s)\frac{\exp(d_0+kp_1)}{2\tanh(h/2)}$.
\end{subsubsection}
%
%
\begin{subsubsection}{Contribution of $x_1$}
If $\g_1$ is open, there is only one summand, which
we have already computed. If $\g_1$ is closed,
we obtain
\begin{align*}
\frac{\pa(p,q,r,s)}{\pa\tau_{\xi,x_1}} & =
\sum_{k\geq 0}\frac{\pa(p,q,r,s)}{\pa\tau_{\ti{\xi}(\ti{x}_1^k)}}+
\sum_{k<0}\frac{\pa(p,q,r,s)}{\pa\tau_{\ti{\xi}(\ti{x}_1^k)}}=\\
& = (p,q,r,s)\left(\sum_{k\geq 0}\frac{\exp(-d_0-kp_1)}{2\tanh(h/2)}
-\sum_{k<0}\frac{\exp(d_0+kp_1)}{2\tanh(h/2)}\right)=\\
& = \frac{(p,q,r,s)}{2\tanh(h/2)}\left(
\exp(-d_0)\sum_{k\geq 0}[\exp(-p_1)]^k-\exp(d_0-p_1)\sum_{j\geq 0}[\exp(-p_1)]^j
\right)= \\
& = \frac{(p,q,r,s)}{2\tanh(h/2)}\,
\frac{\exp(-d_0)-\exp(d_0-p_1)}{1-\exp(-p_1)}
\end{align*}
Multiplying and dividing by $\exp(p_1/2)$, we get
\begin{equation}
\frac{\pa(p,q,r,s)}{\pa\tau_{\xi,x_1}}=
\frac{(p,q,r,s)}{2\tanh(h/2)}\,
\frac{\sinh(p_1/2-d(y_1,x_1))}{\sinh(p_1/2)}
\end{equation}
because $d_0=d(\ti{y}_1,\ti{x}_1^0)=d(y_1,x_1)$.
\end{subsubsection}
%
%
\begin{subsubsection}{Contribution of $x_2$}
Because of the symmetry between $\g_1$ and $\g_2$,
we can apply the same argument above to every
point $x_2\in\xi\cap\g_2$.
\end{subsubsection}
%
%
\begin{proof}[End of the proof of Theorem~\ref{thm:derivative}]
Differentiating the relation $\log|(p,q,r,s)|=\log\sinh^2(h/2)$
we get
\[
\frac{dh}{d(p,q,r,s)}=\frac{\tanh(h/2)}{(p,q,r,s)}
\]
Using the above computations and the chain rule
\[
\frac{\pa h}{\pa \tau_\xi}=\frac{dh}{d(p,q,r,s)}\,
\frac{\pa(p,q,r,s)}{\pa\tau_\xi}
\]
we get the result.
\end{proof}
\end{subsection}
%
%
\begin{subsection}{The general case}\label{subsec:general}
It turns out that the result can be extended to the case
when the $\nu$'s are not necessarily right angles and $\xi$
may intersect $\delta$.

A point of intersection $z\in \xi\cap\delta$ is {\it near $\g_i$}
if $\exists x_i\in \xi\cap\g_i$ such that
$[z,y_i]\subset\delta $, $[y_i,x_i]\subset \g_i$ and $[x_i,z]\subset\xi$
are the sides of a geodesic triangle $T_z$ (locally)
embedded in $R$.
We say that $z$ is {\it distant} if it is not near $\g_1$ or $\g_2$.

Suppose that $\g_i$ is a closed geodesic,
$z\in\xi\cap\delta$ is near $\g_i$
with $T_z=(x_i\,y_i\,z)$.
Travelling from $x_i$
in the direction of $z$ (along a side of $T$), consider the maximum
number $r$ of intersections $z=z_1,z_2,\dots,z_r\in\xi\cap\delta$
such that the loop obtained
as a union of the two arcs $[z_j,z_{j+1}]\subset\xi$
and $[z_{j+1},z_j]\subset\delta$ is homotopic
to $\g_i$ for all $j=1,\dots,r-1$.
If $x_i$ does not belong to any $T_z$, then we set
$r=0$.
\begin{remark}
The exceptional case $x_i=y_i$ must always be treated:
\begin{itemize}
\item[-]
as if $x_i$ comes just after $y_i$ according to
the orientation of $\g_i$,
in case $\nu_{x_i}>\pi/2$;
\item[-]
as if $x_i$ comes just before $y_i$,
in case $\nu_{x_i}<\pi/2$.
\end{itemize}
\end{remark}

Define $r(x_i)=r$ if $(y_i\,x_i\,z)$ is a positively oriented
triangle and $r(x_i)=-r$ if $(y_i\,x_i\,z)$ is negatively oriented.
\begin{center}
\begin{figurehere}
\psfrag{de}{$\tilde{\d}$}
\psfrag{x1}{$\tilde{x}_1$}
\psfrag{x2}{$\tilde{x}_2$}
\psfrag{xi}{$\tilde{\xi}$}
\psfrag{y1}{$\tilde{y}_1$}
\psfrag{y2}{$\tilde{y}_2$}
\psfrag{g1}{$\tilde{\g}_1$}
\psfrag{g2}{$\tilde{\g}_2$}
\psfrag{p}{$p$}
\psfrag{q}{$q$}
\psfrag{s}{$s$}
\psfrag{r}{$r$}
\psfrag{T}{$\tilde{T}_z$}
\psfrag{z}{$\tilde{z}$}
\includegraphics[width=0.5\textwidth]{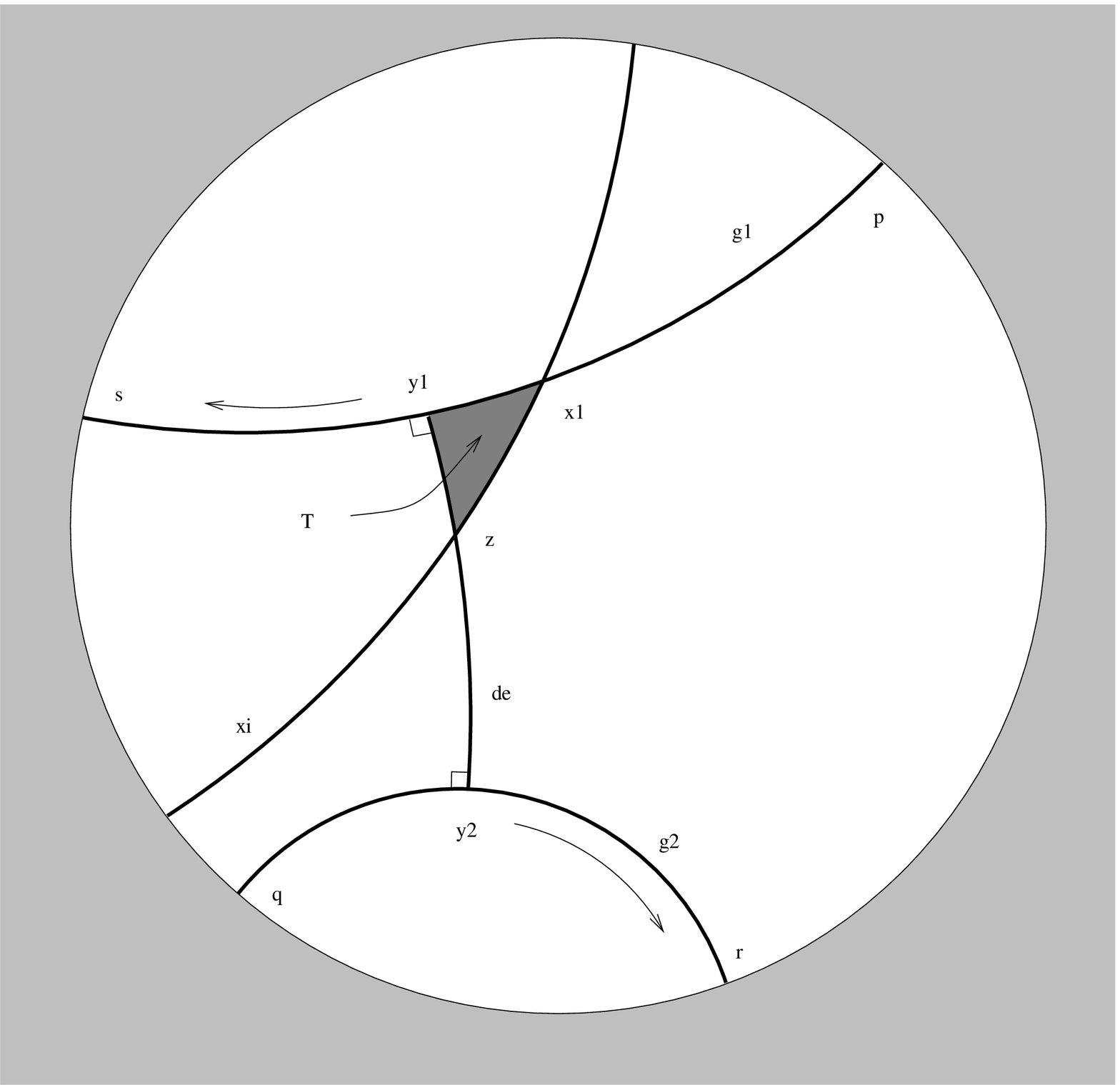}
\caption{Example of lifting of $T_z$ to the universal cover of $R$.}
\end{figurehere}
\end{center}
We will write $[x_i,x_{3-i}]\sim \delta$ (with $i\in\{1,2\}$) if
an oriented segment $[x_i,x_{3-i}]\subset \xi$,
running between $x_i\in\xi\cap\g_i$ and $x_{3-i}\in\xi\cap\g_{3-i}$,
is homotopic to $\delta$ through a homotopy that keeps
the starting point of the segment on $\g_i$ and the endpoint on $\g_{3-i}$.
\begin{center}
\begin{figurehere}
\psfrag{de}{$\d$}
\psfrag{x1}{$ x_1$}
\psfrag{x2}{$x_2$}
\psfrag{xi}{$\xi$}
\psfrag{y1}{$ y_1$}
\psfrag{y2}{$ y_2$}
\psfrag{g1}{$\g_1$}
\psfrag{g2}{$\g_2$}
\psfrag{t1}{$t_{x_1}$}
\psfrag{t2}{$t_{x_2}$}
\includegraphics[width=0.5\textwidth]{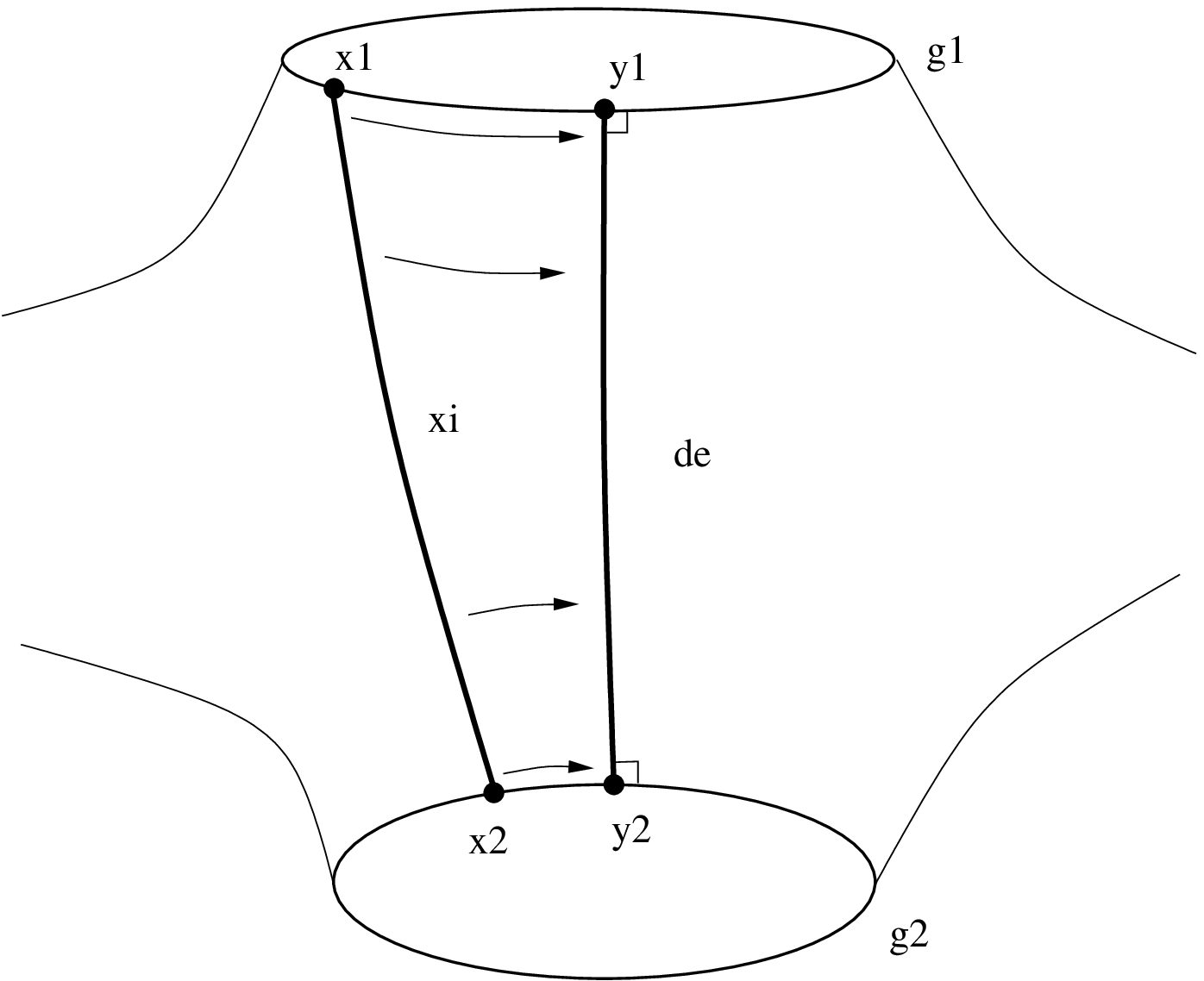}
\caption{Example of $[x_1,x_2]$ homotopic to $\delta$.}
\end{figurehere}
\end{center}
\begin{remark}[On the definition of $d(y_i,x_i)$]
If $\g_i$ is an open geodesic, then the distance (with sign)
$d(y_i,x_i)$ between $y_i\in\delta\cap\g_i$ and $x_i\in\xi\cap\g_i$
is clearly well-defined.

If $\g_i$ is closed, then
\begin{itemize}
\item[-]
if $[x_i,x_{3-i}]\cong\delta$ for some $x_{3-i}\in\g_{3-i}$,
then $d(y_i,x_i)$ is the distance (with sign) between $y_i$
and $x_i$ along the path described by the homotopy that deforms
$[x_i,x_{3-i}]$ to $\delta$;
\item[-]
otherwise, we set $d(y_i,x_i)\in [0,p_i)$.
\end{itemize}
\end{remark}
\begin{theorem}\label{thm:derivative2}
If $\xi$ is any simple closed geodesic, then
\[
\frac{\pa}{\pa\tau_{\xi}}(h)=c_1+c_2+c_0
\]
and we have set
\[
c_0 =\sum_{\substack{z\in\xi\cap\delta \\ \text{distant}}} \cos \a(z)
\]
where $\a(z)$ is the smallest angle one has to rotate the arc of
geodesic $\os{\frown}{z\,y_1}$ (starting at $z$) clockwise
in order to lie on $\xi$;
\[
c_i=\sum_{x_i \in\xi\cap\delta} c_i(x_i)
\]
for $i\in\{1,2\}$ and
\[
c_i(x_i) =
\begin{cases}
\text{$\g_i$ open}
\begin{cases}
\dis\tanh(h/2)\cos(\nu_{x_i}) & \text{if $[x_i,x_{3-i}]\sim\delta$} \\
\phantom{x}\\
\dis \frac{\e\s}{2}\exp\left[-\e\s\, d(y_i,x_i)\right]\sin(\nu_{x_i}) &
\text{otherwise}
\end{cases}\\
\phantom{x}\\
\text{$\g_i$ closed}
\begin{cases}
\dis \frac{\sinh(-d(y_i,x_i))}{2[\exp(p_i)-1]}\sin(\nu_{x_i})
+\tanh(h/2)\cos(\nu_{x_i})
& \text{if $[x_i,x_{3-i}]\sim\delta$}\\
\phantom{x}\\
\dis\frac{\sinh[p_i/2-d(y_i,x_i)-r(x_i)p_i]}{2\sinh(p_i/2)}\sin(\nu_{x_i}) &
\text{otherwise}
\end{cases}
\end{cases}
\]
with $\s=\mathrm{sgn}(d(y_i,x_i))$, $\e=-1$ if $r(x_i)\neq 0$
and $\e=1$ if $r(x_i)=0$.
%
\end{theorem}
The formula above must be compared with Theorem~3.4 in
\cite{wolpert:symplectic}.

The summand $c_0$ comes from distant 
intersections
and is treated in Section~\ref{subsubsec:essential}.

To examine $c_1$ (we will deal similarly with $c_2$),
as before, pick $x_1\in\g_1\cap\xi$ and consider
each branch of $\g_1$ separately,
if $\g_1$ self-intersects at $x_1$.

If $\g_1$ is open, then
the unique lift of $x_1$ along
$\tilde{\g}_1$ will be called $\ti{x}^0_1$
if the lift of $\xi$ through it
separates $s$ from $r$, and $\ti{x}^{-1}_1$
otherwise.

If $\g_1$ is closed, then
let $\ti{x}^0_1\in\ti{\g}_1$ the lift of $x_1$
that separates $s$ from $r$ and which is farthest
from $s$ (in the Euclidean metric of the disc).
Similarly, if $\g_2$ is closed, then $\ti{x}^0_2$
is the lift of $x_2$ that separates $r$ from $s$
and which is farthest from $r$.

If $\g_i$ is closed, consider it as a loop based at $x_i$ and define
$\ti{x}^k_i$ to be the endpoint
of the lift of $(\g_i)^k$ that starts at $\ti{x}^0_i$ for
every $k\in\mathbb{Z}$.
In this case, the distance with sign
$d(\ti{y}_i,\ti{x}^k_i)$ (that is, the length of the portion
of $\ti{\g}_i$ running from $\ti{y}_i$ to $\ti{x}^k_i$
in the positive direction) is exactly $d(y_i,x_i)+(r(x_i)+k)p_i$.
As before, $\ti{\xi}(\ti{x}^k_i)$ is the lift of $\xi$ that
passes through $\ti{x}^k_i$.

The contribution of $x_i\in\xi\cap\g_i$ to
$\pa(p,q,r,s)/\pa\tau_\xi$ is
\[
\frac{\pa(p,q,r,s)}{\pa\tau_{\xi,x_i}}:=
\begin{cases}
\dis\frac{\pa(p,q,r,s)}{\pa\tau_{\ti{\xi}(\ti{x}_i^0)}}
& \text{if $\g_i$ is open and $\ti{\xi}(\ti{x}_i^0)$ separates $s$ and $r$} \\
\dis\frac{\pa(p,q,r,s)}{\pa\tau_{\ti{\xi}(\ti{x}_i^{-1})}}
& \text{if $\g_i$ is open and $\ti{\xi}(\ti{x}_i^{-1})$ does not separate
$s$ and $r$} \\
\dis\sum_{k\in\Z}\frac{\pa(p,q,r,s)}{\pa\tau_{\ti{\xi}(\ti{x}_i^k)}}
& \text{if $\g_i$ is closed}
\end{cases}
\]
and $\dis c_i(x_i)=\frac{\tanh(h/2)}{(p,q,r,s)}\,
\frac{\pa (p,q,r,s)}{\pa\tau_{\xi,x_i}}$.
%
%
\begin{subsubsection}{Case of $\ti{\xi}(\ti{x}^k_1)$
separating $s$ from $\{p,q,r\}$}

By definition, $k\geq 0$.
As in the case with right angles, we have
\[
\pa_{\os{\frown}{s_1 s_2}}(p,q,r,s)=(p,q,r,s)[(q,s_1,s_2,s)-(p,s_1,s_2,s)]
=\frac{(p,q,r,s)}{2}\left(\cos\nu-\cos\th\right)
\]
where $\nu=\nu_{x_1}$.

\begin{center}
\begin{figurehere}
\psfrag{de}{$\tilde{\d}$}
\psfrag{p}{$p$}
\psfrag{q}{$q$}
\psfrag{r}{$r$}
\psfrag{s}{$s$}
\psfrag{n}{$n$}
\psfrag{la}{$\l$}
\psfrag{m}{$m$}
\psfrag{th}{$\th$}
\psfrag{s1}{$ s_1$}
\psfrag{s2}{$ s_2$}
\psfrag{x}{$ \ti{x}_1^k$}
\psfrag{y1}{$ \ti{y}_1$}
\psfrag{y2}{$ \ti{y}_2$}
\psfrag{g1}{$ \ti{\g}_1$}
\psfrag{g2}{$ \ti{\g}_2$}
\psfrag{b}{$\beta$}
\psfrag{nu}{$\nu$}
\includegraphics[width=0.6\textwidth]{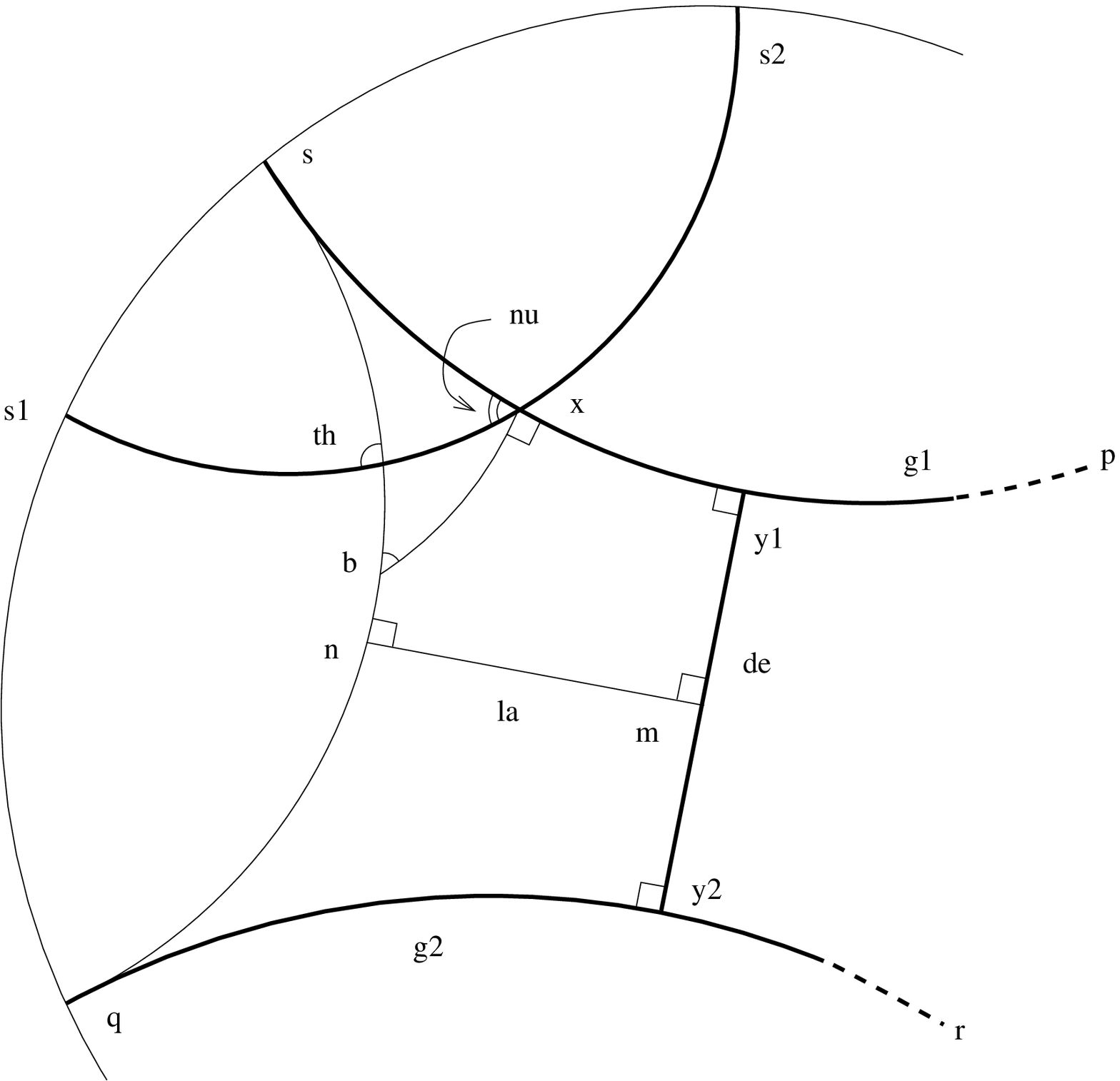}
\end{figurehere}
\end{center}

As before, the picture does not exhaust all the possible cases,
but the formula we will find will hold in all cases, because
of the analyticity mentioned above.

From the previous computations, we know that
$\dis \cos\b=\frac{\exp(-d_k)}{\tanh(h/2)}$.

Call $e$ the length of the segment from $\ti{x}_1^k$
to the vertex of $\b$ and $f$ the length of the segment
from $\ti{x}_1^k$ to the vertex of $\th$.
\[
\cosh e=\frac{\cos\b \cos(\pi/2)+\cos 0}{\sin\b\sin(\pi/2)}=
\frac{1}{\sin\b}
\implies \tanh e=\cos\b
\]
and $\dis \cosh f=\frac{-\cos\nu \cos\th+1}{\sin \nu\sin\th}$
because of part (b) of Lemma~\ref{lemma:triangles} applied
to the triangles $(\b\,\ti{x}_1^k\,s)$ and
$(\th\,\ti{x}_1^k,\,s)$. We also have
\[
\frac{\sin\th}{\sinh e}=\frac{\sin\b}{\sinh f}
\implies
\sin\th \sinh f=\cos\b
\]
because of part (a) of Lemma~\ref{lemma:triangles} applied
to $(\ti{x}_1^k\,\th\,\b)$. From $\sinh^2 f\,\sin^2\th=\cos^2 \b$
we obtain
\[
\frac{\exp(-2d_k)}{\tanh^2(h/2)}=
\frac{(1-\cos\nu\cos\th)^2-\sin^2\nu\sin^2\th}{\sin^2\nu}
\]
Simplifying the expression, we get
\[
(\cos\nu-\cos\th)^2=\frac{\exp(-2d_k)\sin^2\nu}{\tanh^2(h/2)}
\]
As $\nu<\th$, we finally obtain
\[
\cos\nu-\cos\th=\frac{\exp(-d_k)\sin\nu}{\tanh(h/2)}
\]
and so
$\dis
\pa_{\os{\frown}{s_1 s_2}}(p,q,r,s)
=\frac{(p,q,r,s)}{2\tanh(h/2)}\exp(-d_k)\sin\nu$.
\end{subsubsection}
%
%
\begin{subsubsection}{Case of $\ti{\xi}(\ti{x}^k_1)$ separating
$p$ from $\{s,q,r\}$}
By definition, $k<0$. Arguing as in the previous case,
\[
\cos\th+\cos\nu = -\frac{\exp(-d_k)\sin\nu}{\tanh(h/2)}
\]
Hence, we obtain
\begin{align*}
\pa_{\os{\frown}{s_1 s_2}}(p,q,r,s) & =(p,q,r,s)[(p,s_1,s_2,r)-(p,s_1,s_2,s)]=\\
& =\frac{(p,q,r,s)}{2}\left(\cos\th +\cos\nu\right)=\\
& =-\frac{(p,q,r,s)}{2\tanh(h/2)}\exp(-d_k)\sin\nu
\end{align*}
\end{subsubsection}
%
%
\begin{subsubsection}{Case of $\ti{\xi}(\ti{x}^0_1)=\ti{\xi}(\ti{x}^0_2)$
separating $\{p,r\}$ from $\{q,s\}$}\label{subsubsec:c_h}
This happens when the segment $[x_1,x_2]\subset\xi$ is homotopic to $\delta$.
By Lemma~\ref{lemma:cr}
\[
\pa_{\os{\frown}{s_1 s_2}}(p,q,r,s)=(p,q,r,s)[(s,s_1,s_2,q)-
(r,s_1,s_2,q)+(q,s_1,s_2,s)-(p,s_1,s_2,s)]
\]
which gives $\dis \pa_{\os{\frown}{s_1 s_2}}(p,q,r,s)=(p,q,r,s)
(\cos\nu_{x_1}+\cos\nu_{x_2})$.
\end{subsubsection}
%
%
\begin{subsubsection}{The terms $c_i(x_i)$ when $\g_i$ is closed}
We argue similarly to the case with right angles.
If $\ti{\xi}(\ti{x}^0_i)$ does not
separate $\{p,r\}$ from $\{q,s\}$, then $d_0=d(y_i,x_i)+r(x_i)p_i$,
where $d(y_i,x_i)\in[0,p_i)$, and $d_k=d_0+k p_i$.
\[
\frac{\pa(p,q,r,s)}{\pa\tau_{\xi,x_i}}=
\frac{(p,q,r,s)}{2\tanh(h/2)}\,
\frac{\sinh[p_i/2-d(y_i,x_i)-r(x_i)p_i]}{\sinh(p_i/2)}\sin(\nu_{x_i})
\]

If $\ti{\xi}(\ti{x}^0_i)$ separates $\{p,r\}$ from $\{q,s\}$,
then
\[
\frac{\pa(p,q,r,s)}{\pa\tau_{\xi,x_i}}=
\frac{(p,q,r,s)}{2\tanh(h/2)}\,
\frac{\sinh[-d(y_i,x_i)]}{\exp(p_i)-1}\sin(\nu_{x_i})
+(p,q,r,s)\cos(\nu_{x_i})
\]
where the right summand is the contribution of $\ti{x}^0_i$.
\end{subsubsection}

\begin{subsubsection}{Contribution of distant
intersections}\label{subsubsec:essential}
Suppose $z\in\xi\cap\delta$ is a distant intersection
with angle $\a(z)=\a$
and the situation looks like in the figure below,
when lifted to the universal cover.

\begin{center}
\begin{figurehere}
\psfrag{de}{$\tilde{\d}$}
\psfrag{p}{$p$}
\psfrag{q}{$q$}
\psfrag{r}{$r$}
\psfrag{s}{$s$}
\psfrag{n}{$n$}
\psfrag{la}{$\l$}
\psfrag{m}{$m$}
\psfrag{th}{$\th_{qs}$}
\psfrag{s1}{$s_1$}
\psfrag{s2}{$s_2$}
\psfrag{t}{$t$}
\psfrag{z}{$ \ti{z}$}
\psfrag{y1}{$ \ti{y}_1$}
\psfrag{y2}{$ \ti{y}_2$}
\psfrag{g1}{$ \ti{\g}_1$}
\psfrag{g2}{$ \ti{\g}_2$}
\psfrag{b}{$\beta$}
\psfrag{a}{$\alpha$}
\includegraphics[width=0.6\textwidth]{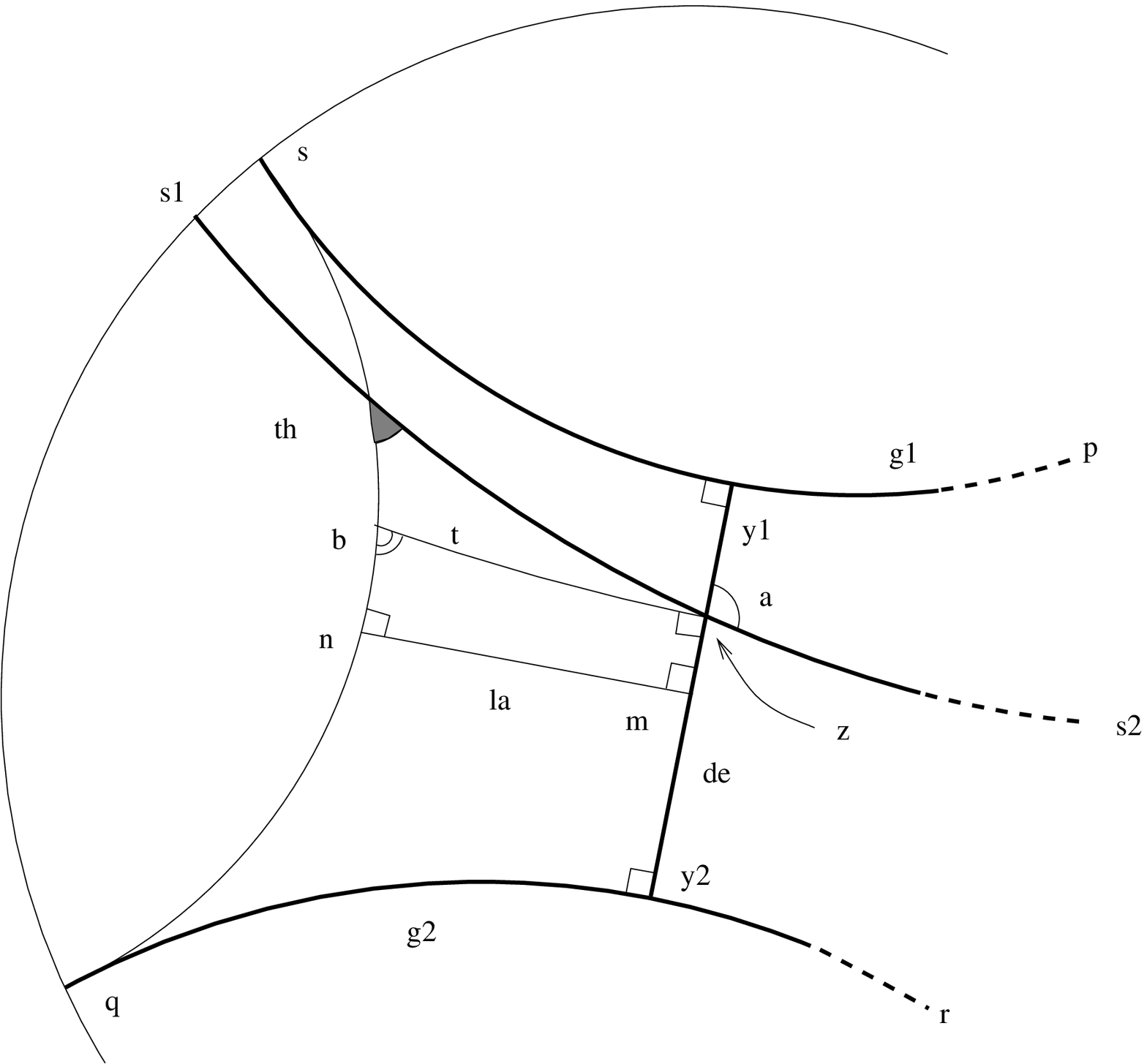}
\label{fig:essential}
\end{figurehere}
\end{center}

Let $e=d(m,\tilde{z})$ be the distance (with sign) between $m$ and $\tilde{z}$,
where $\tilde{\delta}$ is oriented in such a way that
$h=d(\tilde{y}_2,\tilde{y}_1)=-d(\tilde{y}_1,\tilde{y}_2)$.
In Figure~\ref{fig:essential}, we have $e>0$.

In this case, Lemma~\ref{lemma:cr} gives us
\[
\pa_{\os{\frown}{s_1 s_2}}(p,q,r,s)=(p,q,r,s)[(r,s_1,s_2,p)+(q,s_1,s_2,s)-1]
\]
which can be rewritten as
\[
\pa_{\os{\frown}{s_1 s_2}}(p,q,r,s)=-\frac{1}{2}(p,q,r,s)[\cos(\th_{qs})+
\cos(\th_{rp})]
\]
where $\th_{qs}$ is the angle shown in the figure above.

To begin, we have
$\quad\sinh \l\, \sinh(h/2)=1\quad$
from the quadrilateral $(n\,m\,\ti{y}_2\,q)$.
Moreover, $(n\, m\,\tilde{z}\,\beta)$ tells us that
\[
\frac{\sinh e}{\sinh(h/2)}=\sinh e\,\sinh\l=\cos\b
\qquad
\text{and}\qquad
\cosh t=\frac{\cosh\l}{\sin\b}
\]
where $t$ is the length of $\os{\frown}{\b\,\tilde{z}}$.
Looking at
the triangle $(\tilde{z}\,\vartheta_{qs}\,\b)$, we have
\[
\cosh t=\frac{\cos(\pi-\b)\cos(\a-\pi/2)+\cos(\th_{qs})}
{\sin(\pi-\b)\sin(\a-\pi/2)}=
\frac{-\cos\b\, \sin\a+\cos(\th_{qs})}{-\sin\b\,\cos\a}\, .
\]

Because $\dis \cosh\l=\sqrt{1+\sinh^2\l}=\sqrt{1+\sinh^{-2}(h/2)}=
\coth(h/2)$,
we conclude that
\[
\cos(\th_{qs})=\cos\b\,\sin\a-\cosh\l\,\cos\a=
\frac{\sinh e\,\sin\a}{\sinh(h/2)}-\frac{\cos\a}{\tanh(h/2)}\, .
\]

Symmetrically, we have $\dis \cos(\th_{rp})=
\frac{\sinh(-e)\,\sin\a}{\sinh(h/2)}-\frac{\cos\a}{\tanh(h/2)}$
and so
\[
\pa_{\os{\frown}{s_1 s_2}}(p,q,r,s)=\frac{(p,q,r,s)}{\tanh(h/2)}\cos\a\, .
\]
\end{subsubsection}
%
%
\end{subsection}
%
%
\end{section}
%
%
%
\begin{section}{The Weil-Petersson Poisson structure}\label{sec:poisson}
%
%
In this section we want to prove the following.

\begin{theorem}\label{thm:wp}
Let $S$ be a compact hyperbolic surface with boundary
components $\Cc$ and no cusps.
If $\ua=\{\a_1,\dots,\a_{6g-6+3n}\}$ is a triangulation of $S$,
then the Weil-Petersson bivector field on $\Teich(S)$
at $[f:S\rar\Si]$ can be written as
\[
\eta_S=\frac{1}{4}\sum_{C\in\Cc} \sum_{\substack{y_i\in f(\a_i\cap C) \\
y_j\in f(\a_j\cap C) }}
\frac{\sinh(p_C/2-d_C(y_i,y_j))}{\sinh(p_C/2)}\,
\frac{\pa}{\pa a_i}\wedge\frac{\pa}{\pa a_j}
\]
where $a_i=\ell_{\a_i}$, $p_C=\ell_C$ and $d_C(y_i,y_j)\in(0,p_C)$ is
the length of geodesic arc running from $y_i$ to $y_j$ along
$f(C)$ in the positive direction.
\end{theorem}
\begin{remark}
The statement of the theorem still holds if we consider
surfaces with boundary not consisting only of cusps, that
is if we work on $\eTeich(S)\setminus\eTeich(S)(0)$.
In this case, when computing the bivector field at the point
$[f:S\rar\Si]$, one must use a triangulation adapted to $\Sigma$
and the sum involves only arcs of $S$ whose image
through $f$ does not meet the cusps of $\Sigma$.
\end{remark}
\begin{proof}
The triangulation $\ua$ of $S$
determines a pair of pants decomposition
$\{\hat\a_1,\dots,\hat\a_{6g-6+3n}\}$
of the double $dS$, where $\hat\a_i$ is the double of $\a_i$.
As usual, let $\i,\i':S\hra dS$ be the two inclusions and
$D:\Teich(S)\rar\Teich(dS)$ the doubling map induced by $\i$.

Suppose the arc $\a_i$ joins the boundary components $C_s$ and $C_t$ of $S$.
Then, the function $a_i:\Teich(dS)\rar\R_+$ that measures
the length of
the shortest path homotopic to $\i(\a_i)$ that joins the closed geodesics
freely homotopic to $\i(C_s)$ and $\i(C_t)$ reduces to the usual
$a_i$, when
restricted to $D(\Teich(S))$. Similarly, we can define $a'_i:
\Teich(dS)\rar\R_+$ as the length of the shortest path homotopic to
$\i'(\a_i)$ that joins the closed geodesics homotopic to $\i'(C_s)$ and
$\i'(C_t)$.

By Wolpert's theorem (see Section~\ref{sec:2vector}),
the Weil-Petersson bivector field on $\Teich(dS)$
can be written as
\[
\eta_{dS}=-\sum_{i=1}^{6g-6+3n}\frac{\pa}{\pa \ell_i}\wedge
\frac{\pa}{\pa\tau_i}
\]
where $\ell_i=\ell_{\hat\a_i}$ and $\tau_i$ is the twist parameter
associated to $\hat\a_i$. It is immediate to realize that
$D^*(d\tau_i)=0$; really, we can fix the conventions about
the twist coordinates in such a way
that $\tau_i\Big|_{D(\Teich(S))}\equiv 0$.

Proposition~\ref{prop:2vector} tells us that
$\dis(\pi_\i)_*\left(\eta_{dS}\Big|_{\Teich(S)}\right)=\eta_S$
is the Weil-Petersson bivector field
on $\Teich(S)$, where $\pi_\i:\Teich(dS)\rar\Teich(S)$ associates
to $[f:dS\rar R]$ the ``half'' of the surface $R$ corresponding
to $f(\i(S))$.

Now, let's consider the following diagram
\[
\xymatrix{
0\ar[r] & T_{d\Si}\Teich(dS) \ar[r]^{\varphi\qquad\qquad\qquad}
\ar[rd]_{\pi_\i} &
T_\Si\Teich(S)\oplus T_\Si\Teich(S)'\, \oplus
\left(\dis\bigoplus_C \dis\R\frac{\pa}{\pa \tau_C} \right)
\ar[d]^{\pi_1} \ar[r]
& \dis\bigoplus_C \dis \R \frac{\pa}{\pa p_C} \ar[r] & 0 \\
&& T_\Si\Teich(S)
}
\]
where $\pi_1$ is the projection onto the first summand.
Clearly, at $d\Si$
\[
\varphi\left(\frac{\pa}{\pa\tau_i}\right)=
\sum_j \frac{\pa a_j}{\pa\tau_i}\,\frac{\pa}{\pa a_j}+
\sum_k \frac{\pa a'_k}{\pa\tau_i}\,\frac{\pa}{\pa a'_k}+
\sum_{C\in\Cc} \frac{\pa\tau_C}{\pa\tau_i}\,\frac{\pa}{\pa\tau_C}
\]
using the bases $\{\pa/\pa a_i\}$ for $T_\Si\Teich(S)$
and $\{\pa/\pa a'_i\}$ for $T_\Si\Teich(S)'$.
Hence
\[
(\pi_\i)_*\frac{\pa}{\pa\tau_i}=
\sum_j \frac{\pa a_j}{\pa\tau_i}\,\frac{\pa}{\pa a_j}
\]
Moreover, $\dis \varphi\left(\frac{\pa}{\pa\ell_i}\right)=
\frac{1}{2}\left(\frac{\pa}{\pa a_i}+\frac{\pa}{\pa a'_i} \right)$
implies that
$\dis(\pi_\i)_*\frac{\pa}{\pa\ell_i}=\frac{1}{2}
\frac{\pa}{\pa a_i}$.
As a consequence, we deduce
\[
\eta_S=-\frac{1}{2}\sum_{i,j}
\frac{\pa a_j}{\pa\tau_i}\,
\frac{\pa}{\pa a_i}
\wedge\frac{\pa}{\pa a_j}
\]

Given an oriented arc $\ora{\a_i}$ on $S$, call
$y(\ora{\a_i})$ the
endpoint of $\ora{\a_i}$ and $C(\ora{\a_i})$ the boundary component
that contains $y(\ora{\a_i})$. At $[f:S\rar\Si]\in\Teich(S)$,
we will denote by $f(y(\ora{\a_i}))$ the endpoint of the geodesic
arc in the class of $f(\ora{\a_i})$.

Given two distinct oriented arcs $\ora{\a_i}$ and $\ora{\a_j}$
that end on the same component $C=C(\ora{\a_i})=C(\ora{\a_j})$,
the distance $d_C(\ora{\a_j},\ora{\a_i})$ at $[f]$ is the length of the path
from $f(y(\ora{\a_j}))$ to $f(y(\ora{\a_i}))$ along $f(C)$ in the positive
direction.

Applying Theorem~\ref{thm:derivative}, we obtain
\[
\frac{\pa a_j}{\pa\tau_i}=
\frac{1}{2}\sum_{\substack{ \ora{\a_j},\ora{\a_i} \\
C=C(\ora{\a_j})=C(\ora{\a_i}) }}
\frac{\sinh(p_C/2-d_C(\ora{\a_j},\ora{\a_i}))}{\sinh(p_C/2)}
\]
where $p_C$ is the length of $f(C)$. Thus, we have
\begin{align*}
\eta_S & =-\frac{1}{2}\sum_i \frac{\pa}{\pa a_i}\wedge\frac{\pa}{\pa\tau_i}=
-\frac{1}{2}\sum_{i,j}\left(\frac{\pa a_j}{\pa \tau_i}\right)
\frac{\pa}{\pa a_i}\wedge\frac{\pa}{\pa a_j}=\\
& = \frac{1}{4}\sum_{\substack{\ora{\a_i},\ora{\a_j} \\
C=C(\ora{\a_i})=C(\ora{\a_j})}}
\frac{\sinh(p_C/2-d_C(\ora{\a_i},\ora{\a_j}))}{\sinh(p_C/2)}\,
\frac{\pa}{\pa a_i}\wedge\frac{\pa}{\pa a_j}
\end{align*}
which can also be rewritten as
\[
\eta_S=\frac{1}{4}\sum_{C\in\Cc} \sum_{\substack{y_i\in f(\a_i\cap C) \\
y_j\in f(\a_j\cap C) }}
\frac{\sinh(p_C/2-d_C(y_i,y_j))}{\sinh(p_C/2)}\,
\frac{\pa}{\pa a_i}\wedge\frac{\pa}{\pa a_j}
\]
\end{proof}
%
%
\begin{subsection}{The case of large boundary lengths}
Let $\ua=\{\a_i\}_{i=1}^N$ be a triangulation of $S$ and suppose that
there is a sequence of points $[f_n:S\lra\Si_n]\in\Teich(S)$
such that $a_i^{(n)}=\ell_{\a_i}(f_n)\rar 0$
(and so $s(\a_i)^{(n)}=\cosh(a_i^{(n)})\rar 1$)
for all $i$ as $n\rar +\infty$.
We want to study the limit of $\eta_S$ at $[f_n]$ as $n\rar +\infty$.

Let $\ora{\a_i},\ora{\a_j},\ora{\a_k}$ be the oriented arcs
that bound a hexagon in $S\setminus\bigcup_t \a_t$.
\begin{center}
\begin{figurehere}
\psfrag{zi}{$z_i$}
\psfrag{yi}{$y_i$}
\psfrag{mi}{$m_i$}
\psfrag{fi}{$f_i$}
\psfrag{zj}{$z_j$}
\psfrag{yj}{$y_j$}
\psfrag{mj}{$m_j$}
\psfrag{fj}{$f_j$}
\psfrag{zk}{$z_k$}
\psfrag{yk}{$y_k$}
\psfrag{mk}{$m_k$}
\psfrag{fk}{$f_k$}
\psfrag{u}{$u$}
\psfrag{gi}{$\g_i$}
\psfrag{gj}{$\g_j$}
\psfrag{gk}{$\g_k$}
\includegraphics[width=0.7\textwidth]{width.eps}
\end{figurehere}
\end{center}
Remark~\ref{rmk:angle} gives
$\dis \cos(\g_i)=\frac{s(\a_j)^2+s(\a_k)^2-s(\a_i)^2}
{2s(\a_j) s(\a_k)}\rar \frac{1}{2}$.

From $\sinh(w_{\ua}(\ora{\a_i}))\sinh(a_i/2)=\cos(\g_i)$,
we also have $a_i \exp(w_i/2)\rar 2$,
where $w_i=w_{\ua}(\a_i)$.

Hence, $\dis w_i\asymp -2\log(a_i/2)$ and
$\dis
\frac{\pa}{\pa a_i}\asymp
-\exp(w_i/2)\frac{\pa}{\pa w_i}$
in the limit.

Let $\ora{\a_i}$ and $\ora{\a_q}$ be arcs whose endpoints
belong to the same boundary component $C$ and
suppose that the positive path along $C$
from the endpoint
of $\ora{\a_i}$ to the endpoint of $\ora{\a_q}$
meets the endpoints of the oriented arcs
$\ora{\a_{i_0}},\ora{\a_{i_1}},\dots,\ora{\a_{i_l}}$
(where $i_0=i$, $i_l=q$ and we use the convention
$\ola{\a_k}=\ora{\a_{-k}}$).
Then,
$\dis d_C(\ora{\a_i},\ora{\a_q})=
\sum_{r=1}^l d_C(\ora{\a_{i_{r-1}}},\ora{\a_{i_r}})$ and
in the limit
$\dis d_C(\ora{\a_i},\ora{\a_q})\asymp
\frac{w_i+w_q}{2}+\sum_{r=1}^{l-1} w_r$.

Also,
$\dis \frac{\sinh(p_C/2-d_C(\ora{\a_i},\ora{\a_q}))}
{\sinh(p_C/2)}\asymp \exp(-d_C(\ora{\a_i},\ora{\a_q}))-
\exp(d_C(\ora{\a_i},\ora{\a_q})-p_C)$.

Let's compute the limit of the contribution of
$(\ora{\a_i},\ora{\a_q})$ to $\eta_S$.

If $\ora{\a_q}$ comes just after $\ora{\a_i}$ (that is,
$\ora{\a_q}=\ola{\a_j}$ in the picture),
then we obtain a contribution
\[
\asymp \frac{1}{4}\left[\exp\left(-\frac{w_i+w_q}{2}\right)-
\exp\left(\frac{w_i+w_q}{2}-p_C\right)\right]
\exp\left(\frac{w_i+w_q}{2} \right)
\frac{\pa}{\pa w_i}\wedge\frac{\pa}{\pa w_q}
\]
which tends to $0$ if $\ora{\a_i}$ and $\ora{\a_q}$ are
the only oriented arcs incident on $C$, and tends to
$\dis \frac{1}{4}\frac{\pa}{\pa w_i}\wedge\frac{\pa}{\pa w_q}$
otherwise.
We get a similar result if $\ora{\a_q}$ is the oriented arc
that comes just before $\ora{\a_i}$ along $C$.

On the contrary, if $\ora{\a_i}$ and $\ora{\a_q}$ are not
adjacent, then we get a contribution
\begin{multline*}
\asymp \frac{1}{4}
\left[\exp\left(-\frac{w_i+w_q}{2}-\sum_{r=1}^{l-1}w_r\right)-
\exp\left(\frac{w_i+w_q}{2}+\sum_{r=1}^{l-1}w_r-p_C\right)\right]\cdot\\
\cdot\exp\left(\frac{w_i+w_q}{2}\right)
\frac{\pa}{\pa w_i}\wedge\frac{\pa}{\pa w_q}
\end{multline*}
whose coefficient
tends to zero, because all the $w_r$'s diverge in the limit.

Let's use the following normalization:
$\dis\ti{w}_i:=\frac{2w_i}{\sum_C p_C}$, so that
$\dis\sum_i \ti{w}_i=1$.
Then, $\dis \frac{\pa}{\pa w_i}=\frac{2}{\sum_C p_C}\frac{\pa}{\pa\ti{w}_i}$
and we obtain the following.

\begin{theorem}
Let $\ua=\{\a_i\}$ is a triangulation of $S$ and suppose that
there is a sequence of points $[f_n:S\lra\Si_n]\in\Teich(S)$
such that $a_i^{(n)}=\ell_{\a_i}(f_n)\rar 0$
for all $i$ as $n\rar +\infty$.
Call $\eta_S^{(n)}=(\eta_S)_{[f_n]}$ and let $\ti{\eta}_S^{(n)}=
\left(\frac{1}{2}\sum_C p_C^{(n)}\right)^2\eta_S^{(n)}$. Then
\[
\lim_{n\rar\infty} \widetilde{\eta}_S^{(n)}=
\frac{1}{2}\sum_{h\in H}
\left(\frac{\pa}{\pa \ti{w}_i}\wedge\frac{\pa}{\pa \ti{w}_j}
+\frac{\pa}{\pa \ti{w}_j}\wedge\frac{\pa}{\pa \ti{w}_k}
+\frac{\pa}{\pa \ti{w}_k}\wedge\frac{\pa}{\pa \ti{w}_i}\right)
\]
where $H$ is the collection of hexagons in $S\setminus\bigcup_t \a_t$
and $(\a_i,\a_j,\a_k)$ is the cyclically ordered triple of arcs
that bound $h\in H$.
\end{theorem}

We can also compute the class of the limit
of $\widetilde{\eta}_S^{(n)}$ in a
different way.

The following observations are due to Mirzakhani \cite{mirzakhani:virasoro}.

The space $\Teich^*(S)$ is the Teichm\"uller space of
surfaces with $m$ boundary components,
together with a marked point $x_C$ at each boundary
circle $C$. If we have points at the boundaries, we can still
define a $\theta_C$ for each $C$, so that the symplectic
Weil-Petersson form can be restored as
$\dis \sum_i d\ell_i\wedge d\tau_i+\sum_C dp_C\wedge d\theta_C$.

The forgetful map $\Teich^*(S)\lra\Teich(S)$
is clearly a principal $\mathbb{T}=(S^1)^m$-bundle
and the function $\mu=\frac{1}{2}\Ll^2:\Teich^*(S)
\lra\R_{\geq 0}^m$ is a moment map for the action of $\mathbb{T}$.
Notice that all values of $\mu$ are regular.
In fact,
the Teichm\"uller space $(\Teich(S)(\up),\omega_{\up})$ is recovered
as the symplectic
reduction $\mu^{-1}(p_1^2/2,\dots,p_m^2/2)/\mathbb{T}$.

Hence,
using the coisotropic embedding theorem
(see \cite{guillemin:momentmaps}, for instance),
Mirzakhani could conclude that
the following cohomological identity holds
\[
[\omega_{\up}]=[\omega_{0}]+
\frac{1}{2}\sum_C p_C^2 \psi_C
\]
where $\psi_C$ is the first Chern class of the circle bundle over $\Teich(S)$
associated to $C$. Call $\ti{\omega}_{\up}$ the class obtained
dividing $\omega_{\up}$ by $\dis\left(\frac{1}{2}\sum_C p_C\right)^2$.

As $\ti{\omega}_{\up}$ is dual to
$\ti{\eta}_{\up}$, we are interested in computing
\[
[\ti{\omega}_{\up}]=
\frac{4[\omega_{0}]+2\sum_C p_C^2\psi_C}{(\sum_C p_C)^2}
\asymp
\frac{1}{2}\sum_C \ti{p}_C^2\psi_C
\]
where $\dis \ti{p}_C=\frac{2p_C}{\sum_i p_i}$.

However, the argument above involves cohomology classes:
we would like to obtain a pointwise statement.

Theorem~\ref{thm:decomposition} gives us a homeomorphism
$\Phi:\eTeich(S)\setminus\eTeich(S)(0)\lra|\Ao(S)|_{\R}$.
The cells $|\ua|^\circ\subset|\Ao(S)|_{\R}$ 
have affine coordinates $\{e_i\}$, where $e_i$ is the weight of
$\a_i\in\ua$, and $\Phi^*(e_i)=w_i$.

Kontsevich \cite{kontsevich:intersection} wrote a piecewise-linear
$2$-form $\Omega$ on $|\Ao(S)|_\R$
representing (the pull-back from $\M(S)$ of)
$\sum_C p_C^2\psi_C$ and a piecewise-linear
bivector field $\beta$, which is the dual of $\Omega/4$.
The expression of $\beta$
on the top-dimensional cells is the following
\[
\beta=
\sum_{h\in H} \left(\frac{\pa}{\pa e_i}\wedge\frac{\pa}{\pa e_j}
+\frac{\pa}{\pa e_j}\wedge\frac{\pa}{\pa e_k}
+\frac{\pa}{\pa e_k}\wedge\frac{\pa}{\pa e_i}\right)
\]
and its normalized version is
\[
\ti{\beta}=\frac{4\beta}{\left(\sum_C p_C\right)^2}=
\sum_{h\in H} \left(\frac{\pa}{\pa \ti{e}_i}\wedge\frac{\pa}{\pa \ti{e}_j}
+\frac{\pa}{\pa \ti{e}_j}\wedge\frac{\pa}{\pa \ti{e}_k}
+\frac{\pa}{\pa \ti{e}_k}\wedge\frac{\pa}{\pa \ti{e}_i}\right)
\]
where $\dis \ti{e}_t=\frac{2e_t}{\sum_C p_C}$ and $p_C$ is the sum
of the weights of the arcs incident on $C$.
By direct comparison of the explicit expressions for
$\ti{\beta}$ and $\ti{\eta}$, we have the following.

\begin{corollary}
As $n\rar\infty$,
the following limits
\[
2\Phi_*\ti{\eta}_{[f_n]}\rar \ti{\beta}
\quad\text{and}\quad
2\Phi_*\ti{\omega}_{[f_n]}\rar \ti{\Omega}
\]
hold pointwise, where $\dis\ti{\Omega}=\frac{4\Omega}
{\left(\sum_C p_C\right)^2}$.
\end{corollary}
\end{subsection}
%
%
\begin{subsection}{The case of small boundary lengths}
Let $S$ be a Riemann surface with boundary components $\mathcal{C}=\{C_1,\dots,C_m\}$
and $\chi(S)<0$. Remember that $\Ll:\eTeich(S)\rar\R_{\geq 0}^m$ is
the boundary length map, so that $\eTeich(S)(0)=\Ll^{-1}(0)$
is the locus of the surfaces
with $m$ cusps.

Penner has computed the pull-back through
the forgetful map $F:\eTeich(S)(0)\times\R_+^m
\lra\eTeich(S)(0)\subset\eTeich(S)$ of the Weil-Petersson form.

\begin{theorem}[\cite{penner:volumes}]
Fix a triangulation $\ua=\{\a_i\}$ of $S$ and let
$\ti{a}_i:\eTeich(S)(0)\times\R_+^m\lra\R_+$
be the reduced length function
$([f:S\rar\Si],\up)\mapsto \ell_{\a_i}^{\up}(f)$.
Then the pull-back $F^*\omega$ of the
Weil-Petersson $2$-form coincides with
\[
\omega_P:=-\frac{1}{2}\sum_{t\in H} \left(
d\ti{a}_i\wedge d\ti{a}_j+d\ti{a}_j\wedge d\ti{a}_k+
d\ti{a}_k\wedge d\ti{a}_i \right)
\]
where $H$ is the set of hexagons in $S\setminus
\bigcup_i \a_i$ and $(\a_i,\a_j,\a_k)$ is the set of cyclically ordered
arcs that bound the hexagon $t$.
\end{theorem}

\begin{remark}
Using the obvious embedding $(\Delta^\circ)^{m-1}\hra\R_+^m$,
we can pull the functions $\ti{a}_i$'s and
$\omega_P$ back on $\eTeich(S)(0)\times(\Delta^\circ)^{m-1}$.
However, when we regard $(\Delta^\circ)^{m-1}$ as $\R_+^m/\R_+$,
the natural coordinates on $\eTeich(S)(0)\times(\Delta^\circ)^{m-1}$
are the differences $(\ti{a}_i-\ti{a}_{i_0})_{i\neq i_0}$
for any fixed $i_0$.
Notice that a different choice of the constant $M>0$ used for the
embedding $(\Delta^\circ)^{m-1}\hra\{\up\in\R_+^m\,|\, p_1+\dots+p_m=M\}\subset\R_+^m$
will just produce a shift $\ti{a}_i\mapsto \ti{a}_i+\log M$.
Thus, the differences $\ti{a}_i-\ti{a}_{j}$, the
$d\ti{a}_i$'s and $\omega_P$ on $\eTeich(S)(0)\times(\Delta^\circ)^{m-1}$
are well-defined.
\end{remark}

For every $([f],\up)\in\eTeich(S)(0)\times\R_+^m$, define 
\[
\eta_{[f],\up}=
\frac{1}{4}\sum_{C\in\Cc} \sum_{\substack{y_i\in f(\a_i\cap C) \\
y_j\in f(\a_j\cap C) }}
\left(1-\frac{2d_C(y_i,y_j)}{p_C}\right)
\frac{\pa}{\pa \ti{a}_i}\wedge\frac{\pa}{\pa \ti{a}_j}
\]
where $d_C(y_i,y_j)\in(0,p_C)$ is the distance along the
horocycle corresponding to $f(C)$.
It descends to a bivector field on $\eTeich(S)(0)\times(\Delta^\circ)^{m-1}$
and describes the extension of $\eta$ over the real blow-up
$\mathrm{Bl}_{\eTeich(S)(0)}\eTeich(S)$ (whose fiber over $\eTeich(S)(0)$
can be identified to $\eTeich(S)(0)\times\Delta^{m-1}$),
because 
$(a_i-a_j)(f_n)\rar (\ti{a}_i-\ti{a}_j)(f)$
for all $i,j$ as $n\rar \infty$.
\begin{proposition}\label{prop:short}
Fix a triangulation $\ua$ of $S$.
For every $([f],\up)\in\eTeich(S)(0)\times\R_+^m$
\[
\omega_P(\eta_{[f],\up}(d\ti{a}))=
d\ti{a}+d\log(p_+)+d\log(p_-)
\]
where $\ti{a}=\ell^{\up}_\a$ and
$\a\in\ua$ joins $C_+$ and $C_-$.
\end{proposition}

Fix a surface with a projective decoration
$([f:S\rar\Si],[\up])\in\eTeich(S)(0)\times(\Delta^\circ)^{m-1}$,
where $\up=(p_1,\dots,p_m)\in\R_+^m$, such that
$\Sp^*_{fin}(\Si,\up)$ is a triangulation, and
call $c_i=f(C_i)$ the $i$-th cusp of $\Si$.
Consider a sequence of points $[f_n:S\lra\Si_n]\in\Teich(S)$
such that $([f_n],[\up^{(n)}])$ converges to $([f],[\up])$
in $\eTeich(S)\times\Delta^{m-1}$ as $n\rar +\infty$,
where $\up^{(n)}=\Ll(f_n)$.

\begin{corollary}
The limit symplectic structure along the leaf $\Ll^{-1}(\up^{(n)})$
at $[f_n]$ dual to $\eta$ converges to $\omega_{P}$ at
$([f],[\up])\in\eTeich(S)(0)\times\Delta^{m-1}$ (as it must be).
\end{corollary}

Notice that the assertion follows from Proposition~\ref{prop:short}
and the fact
that the symplectic leaves of $\eTeich(S)(0)\times\R_+^m$
are defined by $dp_1=\dots=dp_m=0$.

\begin{notation}
Let $([f:S\rar\Si],\up)$ be a decorated
hyperbolic surface and let $\Si_{\up}$ be the associated
truncated surface.
For every oriented arc $\ora{\a_i}$ of $\ua$ starting
at the boundary component $C$, call
$e(\ora{\a_i})$ the sum of the lengths of the two horocyclic
arcs running around $f(C)$ from the starting point of $f(\ora{\a_i})
\cap\Si_{\up}$
to the previous and the following arc.
Given a portion $\th_{x}$ of the oriented component $C$
running from $x$ to $x'$, where $x,x'$ are consecutive
points in $\dis P_C:=C\cap\left(\bigcup \a_i\right)$,
then the arc {\it opposed to $\th_{x}$} is
the arc $\th_x^{op}\in\ua$ facing $\th_{x}$ in the truncated triangle
that contains $\th_{x}$. Denote by $f(P_C)$ the corresponding
points of $\pa\Si_{\up}\cap f(\bigcup\a_i)$.
\end{notation}
\begin{proof}[Proof of Proposition~\ref{prop:short}]
Pick a truncated triangle $t$ of $S\setminus\bigcup_{\a\in\ua}\a$
and let $\a_i,\a_j,\a_k\in\ua$
be the (cyclically ordered) arcs that bound $t$.
Then the length of the horocyclic arc between $f(\a_j)$
and $f(\a_k)$ is $\dis 2h_{t,i}=\frac{2\l_i}{\l_j\l_k}$.
This implies that $\dis 2\frac{\pa h_{t,i}}{\pa a_i}=h_{t,i}$,
whereas $\dis 2\frac{\pa h_{t,i}}{\pa a_j}=-h_{t,i}$
and $\dis 2\frac{\pa h_{t,i}}{\pa a_k}=-h_{t,i}$.
Because $p_C$ is the sum of all the horocyclic arcs
around $f(C)$ running between consecutive points
of $f(P_C)$,
we easily get
\[
dp_C=-\frac{1}{2}\sum_{\substack{\ora{\a_i}\ \text{out}\\\text{from $C$}}}
e(\ora{\a_i})d\ti{a}_i+\frac{1}{2}\sum_{x\in P_C} h_x dh_x^{op}
\]
where $h_x$ (resp. $h_x^{op}$) is the length of $\th_x$
(resp. $\th_x^{op}$).

Let $\ora{\a}$ be an orientation of the arc
$\a\in\ua$ and let $C_+=C(\ola{\a})$
be the ``source'' of $\ora{\a}$ and call $p_+=\ell_{C_+}$
and $x_0$ the starting point of $\ora{\a}$.
Define similarly $C_-$, $p_-$ and $y_0$.
\begin{center}
\begin{figurehere}
\psfrag{x0}{$x_0$}
\psfrag{x1}{$x_1$}
\psfrag{x2}{$x_2$}
\psfrag{xk}{$x_k$}
\psfrag{xk-1}{$x_{k-1}$}
\psfrag{y0}{$y_0$}
\psfrag{y1}{$y_1$}
\psfrag{y2}{$y_2$}
\psfrag{yl}{$y_l$}
\psfrag{yl-1}{$y_{l-1}$}
\psfrag{a}{$\ora{\a}$}
\psfrag{b1}{$\ora{\b_1}$}
\psfrag{b2}{$\ora{\b_2}$}
\psfrag{bk}{$\ora{\b_k}$}
\psfrag{bk-1}{$\ora{\b_{k-1}}$}
\psfrag{g1}{$\ora{\g_1}$}
\psfrag{g2}{$\ora{\g_2}$}
\psfrag{gl}{$\ora{\g_l}=\th_{x_0}^{op}$}
\psfrag{gl-1}{$\ora{\g_{l-1}}$}
\psfrag{hx0}{$\th_{x_0}$}
\psfrag{C+}{$C_+$}
\psfrag{C-}{$C_-$}
\includegraphics[width=0.7\textwidth]{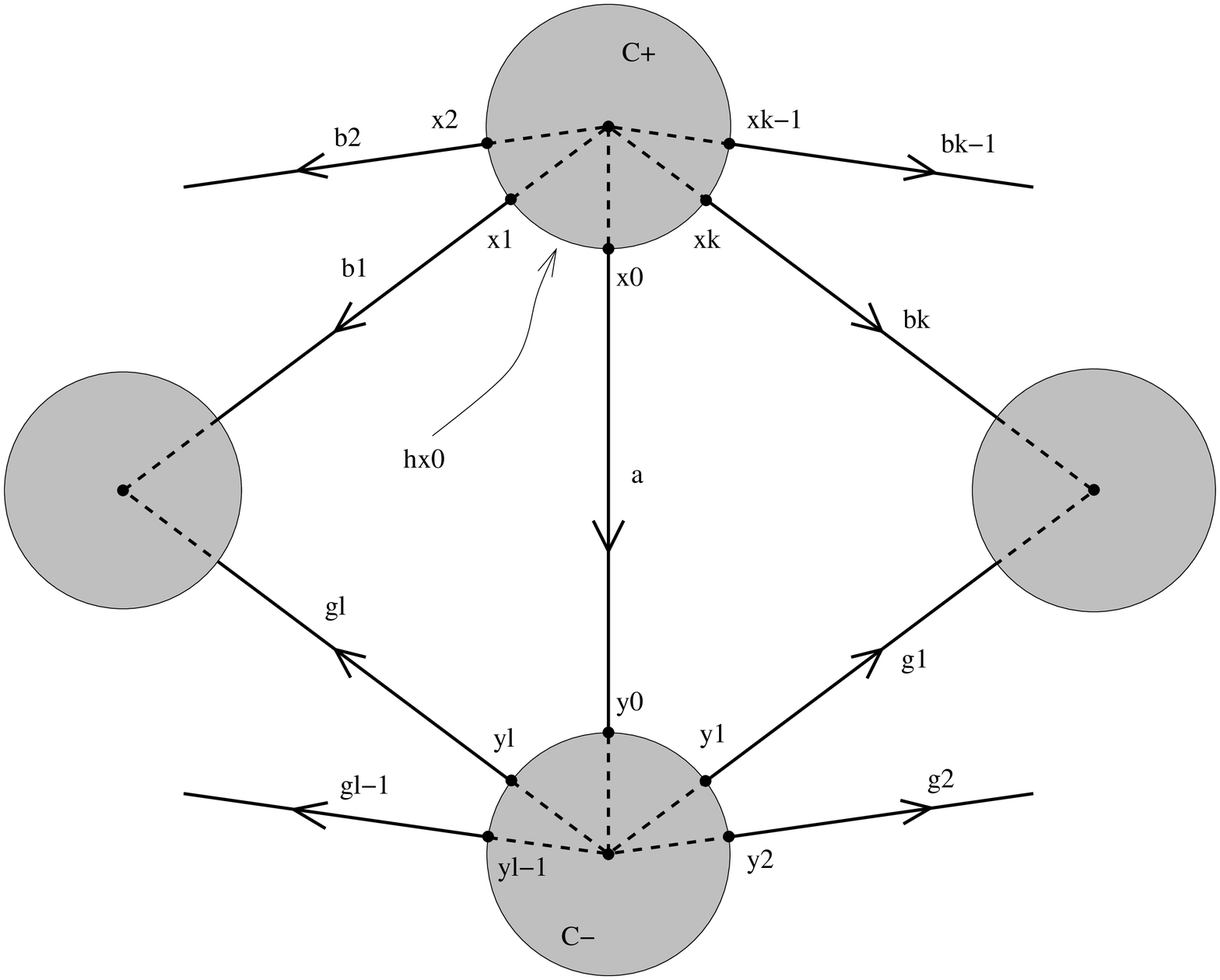}
\end{figurehere}
\end{center}
Starting from $\ora{\a}$ and moving along $C_+$ in the positive
direction, call $\ora{\b_1},\ora{\b_2},\dots,\ora{\b_k}$ the
(ordered) arcs outgoing from $C_+$ and let $x_i$ be
the starting point of $\ora{\b_i}$.
Similarly, call $\ora{\g_1},\dots,\ora{\g_l}$ the arcs outgoing
from $C_-$ and let $y_j$ be their starting point.
Denote by $\ti{b}_i$ the length of $\Si_{\up}\cap f(\b_i)$ and
by $\ti{c}_j$ the length of $\Si_{\up}\cap f(\g_j)$.

In analyzing $\omega\circ\eta(d\ti{a})$, we get four
different contributions:
the contribution to $d\ti{a}$; the contribution to $d\ti{b}_1$
(and similarly to $d\ti{b}_k$, $d\ti{c}_1$, $d\ti{c}_l$);
the contribution to $d\ti{b}_i$ for $i\neq 1,k$
(and similarly to $d\ti{c}_j$ for $j\neq 1,l$);
the contribution to $dh_{x_i}^{op}$ for $i\neq 0,k$
(and similarly to $dh_{y_j}^{op}$ for $j\neq 0,l$).
The other contributions are immediately seen to vanish.

A direct computation shows that
\begin{align*}
\omega\circ\eta(d\ti{a}) & =-\frac{d\ti{a}}{4}\Big[
-\left(1-\frac{2d_{C_+}(x_k,x_0)}{p_+}\right)
-\left(1-\frac{2d_{C_+}(x_0,x_1)}{p_+}\right)+\\
& \qquad\qquad -\left(1-\frac{2d_{C_-}(y_0,y_1)}{p_-}\right)
-\left(1-\frac{2d_{C_-}(y_l,y_0)}{p_-}\right)
\Big]=\\
& = d\ti{a}-\frac{d\ti{a}}{2}\left(
\frac{e(\ora{\a})}{p_+}+\frac{e(\ola{\a})}{p_-}
\right)
\end{align*}
which is exactly the contribution to $d\ti{a}$
of $d\ti{a}+d\log(p_+)+d\log(p_-)$.

Similar computations can be carried over in the
other three cases.
\end{proof}
\end{subsection}
%
%
\end{section}
%
%
\appendix
\begin{section}{Some formulae from hyperbolic
trigonometry}\label{app:hyperbolic}
The following results of elementary hyperbolic trigonometry
are frequently used throughout the paper.
Proofs can be found on \cite{ratcliffe:hyperbolic}.

The first lemma is the statement of the hyperbolic
laws of sines and cosines.

\begin{lemma}\label{lemma:triangles}
Let $A,B,C$ be the vertices of a hyperbolic triangle 
with angles $\a,\b,\g$ (resp. at $A,B,C$).
\begin{itemize}
\item[(a)] {\it (sine law)}
\[
\frac{\sin\a}{\sinh(BC)}=\frac{\sin\b}{\sinh(AC)}=\frac{\sin\g}{\sinh(AB)}
\]
\item[(b)] {\it (cosine law)}
\begin{align*}
\cosh(AB) & =\frac{\cos\a \cos\b+\cos\g}{\sin\a \sin\b}\\
\cos(\a)  & =\frac{\cosh(AB)\cosh(AC)-\cosh(BC)}{\sinh(AB)\sinh(AC)}
\end{align*}
\end{itemize}
\end{lemma}

The following lemma is about quadrilaterals with at least
two right angles.

\begin{lemma}\label{lemma:quadrilateral}
Let $A,B,C,D$ be the vertices of a hyperbolic quadrilateral.
\begin{itemize}
\item[(a)]
If the angles at $A,B,C$ are right, then
\[
\sinh(AB)\cdot\sinh(BC)=\cos(\g)
\]
where $\g$ is the angle at $D$.
\item[(b)]
If the angles at $C$ and $D$ are right, then
\[
\cosh(AB)=\frac{\cos(\alpha)\cos(\beta)+\cosh(CD)}{\sin(\alpha)\sin(\beta)}
\]
where $\alpha$ is the angle at $A$ and $\beta$ is the angle at $B$.
\end{itemize}
\end{lemma}

The next lemma is about pentagons with four right angles.

\begin{lemma}\label{lemma:pentagon}
Let $A,B,C,D,E$ be the vertices of a
hyperbolic pentagon with four right angles
at $A,B,C,D$. Then
\[
\cosh(BC)=\frac{\cosh(AB)\cdot\cosh(CD)+\cos(\g)}
{\sinh(AB)\cdot\sinh(CD)}
\]
where $\g$ is the angle at $E$ (which is thus
opposed to $BC$).
\end{lemma}

The last lemma deals with the well-known
case of hexagons with six right angles.

\begin{lemma}\label{lemma:hexagon}
Let $A,B,C,D,E,F$ be the vertices of a hyperbolic
hexagon with six right angles. Then
\[
\cosh(BC)=\frac{\cosh(AB)\cdot\cosh(CD)+\cosh(EF)}
{\sinh(AB)\cdot\sinh(CD)}.
\]
\end{lemma}
\end{section}
%
%
%
%
\bibliographystyle{amsalpha}
\bibliography{biblio-wp}
\end{document}